\documentclass[11pt]{article}

\usepackage{graphicx}
\RequirePackage[OT1]{fontenc}
\RequirePackage{amsthm,amsmath,amssymb}

\usepackage{algorithmic}
\usepackage{algorithm}
\usepackage{array}
\usepackage[caption=false,font=normalsize,labelfont=sf,textfont=sf]{subfig}
\usepackage{textcomp}
\usepackage{stfloats}
\usepackage{url}
\usepackage{verbatim}
\usepackage{graphicx}
\usepackage{cite}

\newtheorem{thm}{Theorem}

\newtheorem{lem}{Lemma}

\newcommand{\var}{{\rm Var}}

\newcommand{\cA}{{\cal A}}

\newcommand{\bX}{{\bf X}}

\newcommand{\bp}{{\bf p}}

\newcommand{\bq}{{\bf q}}
\newcommand{\wht}{\widehat}

\newcommand{\wbox}{\sqcap\llap{$\sqcup$}}

\newcommand{\bmu}{\boldsymbol{\mu}}
\newcommand{\btau}{\boldsymbol{\tau}}

\begin{document}

\title{Likelihood Scores for Sparse Signal and Change-Point Detection}


\author{Shouri Hu\thanks{School of Mathematical Sciences, University of Electronic Science and Technology of China},
Jingyan Huang\thanks{Department of Statistics and Data Science,
National University of Singapore}, 
Hao Chen\thanks{Department of Statistics,
University of California at Davis}
and Hock Peng Chan\footnotemark[2]}


\maketitle


\begin{abstract}
We consider here the identification of change-points on large-scale data streams.
The objective is to find the most efficient way of combining information across data stream
so that detection is possible under the smallest detectable change magnitude.
The challenge comes from the sparsity of change-points when only a small fraction of data streams undergo change
at any point in time.
The most successful approach to the sparsity issue so far has been the application of hard thresholding such that only
local scores from data streams exhibiting significant changes are considered and added.
However the identification of an optimal threshold is a difficult one.
In particular it is unlikely that the same threshold is optimal for different levels of sparsity.
We propose here a sparse likelihood score for identifying a sparse signal.
The score is a likelihood ratio for testing between the null hypothesis of no change against an alternative hypothesis
in which the change-points or signals are barely detectable.
By the Neyman-Pearson Lemma this score has maximum detection power at the given alternative.
The outcome is that we have a scoring of data streams that is successful in detecting at the boundary of the 
detectable region of signals and change-points.
The likelihood score can be seen as a soft thresholding approach to sparse signal and change-point detection in which
local scores that indicate small changes are down-weighted much more than local scores indicating large changes.
We are able to show second-order optimality of the sparsity likelihood score in the sense of achieving successful 
detection at the minimum detectable order of change magnitude as well as at the minimum detection asymptotic constant
with respect this order of change. 
\end{abstract}


\section{Introduction}

Consider a large number $N$ of data streams containing change-points.
We consider the situation in which all data up to a given time is available for analysis, so
each data stream is an observed sequence of length $T$.
At each change-point one or more of the sequences undergo distribution change.
The objective is to identify these change-points and the sequences undergoing distribution change.
Of interest here is the identification of these change-points when there is sparsity,
that is when the number of sequences undergoing change is small compared to $N$.
More specifically we want to know the minimum magnitude of change for which the distribution change can be detected
under sparsity.
And secondly we want to be have an algorithm that is able to detect,
with high probability, change-points under the minimum detectable
change.
See Niu, Hao and Zhang \cite{NHZ16} and Wang and Samworth \cite{WS18} for applications to 
engineering, genomics and finance.

A typical strategy to deal with sparsity is to subject local scores to thresholding or penalization before 
summing them up across sequence.
Algorithms employing this strategy include the Sparsified Binary Segmentation (SBS) (Cho and Fryzlewicz \cite{CF15}), 
the double CUSUM (DC) (Cho \cite{Cho16}),
the Informative Sparse Projection (INSPECT) (Wang and Samworth \cite{WS18}) and the scan algorithm of Enikeeva and Harachaoui 
\cite{EH19}.
The strategy was also employed by Mei \cite{Mei10}, Xie and Siegmund \cite{XS13}
and Wang and Mei \cite{WM15} in sequential change-point detection on multiple sequences,
and Zhang, Siegmund, Ji and Li \cite{ZSJL10} to detect distribution deviations from known baselines on multiple sequences. 
Thresholding and penalization suppress noise by removing small and moderate scores,
mostly from the majority of sequences without change,
thus enhancing the signals from the sparse sequences with changes.
It is however unlikely that we are able to specify a threshold or penalization parameter that is optimal at all levels of sparsity.

The higher-criticism (HC) test statistic, 
proposed by Tukey \cite{Tuk76} to check for significantly large number of small p-values,
uses multiple thresholds for sparse mixture detection.
The number of p-values below a threshold is transformed to a
higher-criticism score and this score is maximized over all thresholds.
The Berk and Jones \cite{BJ79} test statistic uses multiple thresholds as well but it applies a different
p-value transformation.
The HC test statistic was shown by Donoho and Jin \cite{DJ04} to be optimal in the detection of a sparse normal mixture.
Cai, Jeng and Jin \cite{CJJ11} extended it to detect intervals in multiple sequences where the means of a sparse fraction of the sequences deviate from a known baseline and showed that the HC test statistic is optimal.
Chan and Walther \cite{CW15} considered sequence length much larger than number of sequences 
with detection boundaries that are more complex. 
They showed that the HC test statistic achieves detection at these boundaries and is optimal in more general settings. 
They also showed that the Berk-Jones test statistic achieves the same optimality. 
The HC and Berk-Jones test statistics have the advantage of not requiring a threshold to be specified in advance.
However at each threshold we are unable to differentiate a p-value just below the threshold and a p-value that is below
the threshold by a large margin.
This loss of information results sub-optimality when identifying the exact location where change has occurred.

Our approach here is to convert the p-values into likelihood scores for testing sparse sequences.
It can be considered to be a soft form of thresholding in which p-values that are close to zero are penalized less than 
p-values that are barely significant.
Unlike the HC and Berk-Jones test statistic only one transformation of the p-values is needed.
It retains the advantage of the HC and Berk-Jones test statistics in not specifying a threshold in advance,
but unlike the HC and Berk-Jones test statistics a smaller p-value always has a larger likelihood score.

Since the likelihood scores are transformations of p-values, the proposed method can be applied to
any type of distribution changes and it can handle data types that vary across sequences.
Our theory however requires a specific distribution family for neat asymptotics and we consider here in particular 
either normal or Poisson data.
We show optimality up to the correct asymptotic constants.
For sparse normal change-points these constants are two-dimensional extensions of those in Ingster \cite{Ing97} 
and Donoho and Jin \cite{DJ04} for sparse normal mixture detection.
These constants have been discussed in the context of sparse normal change-point detection assuming a known baseline  
in Chan and Walther \cite{CW15} and Chan \cite{Cha17}. 
For sparse Poisson change-points the constants are new and different from sparse normal constants.

The optimality of multiple sequence identification of change-points up to the correct constant is new.
Previous works on optimality for normal data are up to the correct order of magnitude though they go beyond
the i.i.d. model,
for example Pilliat, Carpentier and Verzelen \cite{PCV20} 
considered sparse change-point detection in time-series with normal errors.
Liu, Guo and Samworth \cite{LGS21} showed optimality up to the best order for normal errors.
However they considered a different problem that assumes at most one change-point.
 
The algorithm we propose here has two steps in the identification of two change-points.
The first detection screening step applies the Screening and Ranking (SaRa) idea of Niu and Zhang \cite{NZ12}.
The second estimation step for more precise location of change-points
uses the CUSUM-like procedure of Wild Binary Segmentation (WBS), cf. Fryzlewicz \cite{Fry14}..
This two-step approach saves computation time because the fast screening step evaluates a large number of segments whereas
the computationally intensive estimation step is only applied when a change-point has been detected during screening.
In contrast for WBS the estimation step is applied on a large number of randomly generated segments.
Unlike in Niu and Zhang \cite{NZ12} we do not apply the BIC criterion of Zhang and Siegmund \cite{ZS07} to determine the number of change-points.
Instead critical values are specified in advance and binary segmentation, cf. Olshen et al. \cite{OVLW03},
is applied to detect the change-points sequentially.

An alternative to binary segmentation is estimating the full set of change-points at one go by applying global optimization
and making use of dynamic programming to manage the computational complexity.
This was employed by the HMM algorithms of Yao \cite{Yao84} and Lai and Xing \cite{LX11},
the multi-scale SMUCE algorithm of Frick, Munk and Sieling \cite{FMS14} 
and the Bayesian Likelihood algorithm of Du, Kao and Kou \cite{DKK16}.
These methods are however designed for single sequence segmentation.
Niu, Hao and Zhang \cite{NHZ16} provides an excellent background of the historical developments.

The outline of this paper is as follows.
In Section II we introduce the sparse likelihood (SL) scores and show that they are optimal in the detection of sparse normal mixtures.
In Section III we extend SL scores to detect change-points in multiple sequences.
In Section IV we show that SL scores are optimal for change-point detection when the observations are normal or Poisson.
In Section V we perform simulation studies on the SL scores.
In the appendices we prove the optimality of SL scores.

\subsection{Notations}

We write $a_n \sim b_n$ to denote $\lim_{n \rightarrow \infty} (a_n/b_n)=1$.  
We write $a_n=o(b_n)$ to denote $\lim_{n \rightarrow \infty} (a_n/b_n)=0$.
We write $a_n \lesssim b_n$ to denote $a_n \leq Cb_n$ for all $n$ for some $C>0$ and $a_n \asymp b_n$ to denote $a_n \lesssim b_n$ and
$b_n \lesssim a_n$.
We write $X_n = O_p(a_n)$ to denote $P(X_n \leq Ca_n) \rightarrow 1$ for some $C>0$. 
Let $\lfloor \cdot \rfloor (\lceil \cdot \rceil)$ denote the greatest (least) integer function.
Let $\phi$ and $\Phi$ denote the density and distribution function respectively of the standard normal.
Let ${\bf 1}$ denote the indicator function.
Let $\emptyset$ denote the empty set and let $\# A$ denote the number of elements in a set $A$.

\section{Sparse mixture detection}

We start with the simpler problem of detecting a sparse mixture, 
with the objective of motivating the sparse likelihood score.

Let $\bp = (p^1, \ldots, p^N)$ be independent p-values of $N$ null hypotheses 
and let $p^{(1)} \leq \cdots \leq p^{(N)}$ be the sorted p-values. 
Tukey (1976) proposed the higher-criticism test statistic 
\begin{equation} \label{HC}
\mbox{HC}(\bp) = \max_{n: Np^{(n)} \leq n} \tfrac{n-N p^{(n)}}{\sqrt{N p^{(n)}(1-p^{(n)})}},
\end{equation}
with ${\rm HC}(\bp)=0$ if $Np^{(n)}>n$ for all $n$,
for the overall test that all null hypotheses are true.

Donoho and Jin \cite{DJ04} showed that the HC test statistic is optimal for detecting a sparse fraction of false null hypotheses.
Consider test scores $Z^n \sim {\rm N}(0,1)$ when the $n$th null hypothesis is true and $Z^n \sim {\rm N}(\mu_N,1)$ for some $\mu_N>0$ when the $n$th null hypothesis is false.
Define 
\begin{equation} \label{rhoZ}
\rho_Z(\beta) = \left\{ \begin{array}{ll} \beta-\tfrac{1}{2} & \mbox{ if } \tfrac{1}{2} < \beta < \tfrac{3}{4}, \cr
(1-\sqrt{1-\beta})^2 & \mbox{ if } \tfrac{3}{4} \leq \beta < 1. \end{array} \right.
\end{equation}
Donoho and Jin \cite{DJ04} showed that on the sparse mixture $(1-\epsilon) {\rm N}(0,1)+\epsilon {\rm N}(\mu_N,1)$
no algorithm is able to achieve,
as $N~\rightarrow~\infty$,
 \begin{equation} \label{zero}
P_0(\mbox{Type I error}) + P_{\mu_N}(\mbox{Type II error}) \rightarrow 0, 
\end{equation}
for testing $H_0$: $\epsilon=0$ versus $H_1$: $\epsilon=N^{-\beta}$,
if $\mu_N = \sqrt{2 \nu \log N}$ for $\nu < \rho_Z(\beta)$.
They also showed that the HC test statistic achieves (\ref{zero}) when $\nu > \rho_Z(\beta)$
and is thus optimal.
Type I error refers to the conclusion of $H_1$ when $H_0$ is true whereas
Type II error refers to the conclusion of $H_0$ when $H_1$ is true.
Ingster (1997, 1998) established the detection lower bound 
showing that (\ref{zero}) cannot be achieved when $\nu < \rho_Z(\beta)$.

Like the HC test statistic,
the Berk and Jones \cite{BJ79} test statistic 
$$\mbox{BJ}(\bp) = \max_{n: Np^{(n)} \leq n} \Big[ n \log \Big( \tfrac{n}{Np^{(n)}} \Big)+(N-n) \log \Big( \tfrac{N-n}{N(1-p^{(n)})} \Big) \Big]
$$
achieves (\ref{zero}) when $\nu > \rho_Z(\beta)$.

We introduce the sparse likelihood scores in Section II.A and show that they achieve (\ref{zero}) in the detection of sparse mixtures, 
when $\nu > \rho_Z(\beta)$,
in Section II.B.

\subsection{Sparse likelihood}

Let $f_1(p) = \tfrac{1}{p(2-\log p)^2} - \tfrac{1}{2}$ and $f_2 (p) = \tfrac{1}{\sqrt{p}}-2$.
For both $i=1$ and 2, 
$\int_0^1 f_i(p) dp=0$ and $f_i(p)$ increases as $p$ decreases.

Define the sparse likelihood score
\begin{eqnarray} \label{PLp}
\ell_N(\bp) & = & \sum_{n=1}^N \ell_N(p^n), \\ \nonumber
\mbox{where } \ell_N(p) & = & \log \Big( 1+ \tfrac{\lambda_1 \log N}{N} f_1(p) + \tfrac{\lambda_2}{\sqrt{N \log N}} f_2(p) \Big), 
\end{eqnarray}
with $\lambda_1 \geq 0$ and $\lambda_2>0$.
For the purpose of sparse mixture detection $f_1$ plays no role and we can consider $\lambda_1=0$.
The importance of $f_1$ comes from detecting very sparse change-points.

The sparse likelihood is the likelihood of  
\begin{eqnarray*}
p^n & \sim_{\rm i.i.d.} & {\rm Uniform}(0,1), \cr
\mbox{ vs } p^n & \sim_{\rm i.i.d.} & f = 1+\tfrac{\lambda_2}{\sqrt{N \log N}} f_2.
\end{eqnarray*}
The distribution function of $f$ is
$$F(p) = p + \tfrac{2 \lambda_2}{\sqrt{N \log N}} (\sqrt{p}-p).
$$
Consider the true distribution of $p^n$ to be some unknown $G \neq {\rm Uniform}(0,1)$.
The fraction of p-values less than some small $p_0$ has mean $G(p_0)$ and standard deviation approximately $\sqrt{G(p_0)/N}$.
Hence to achieve decent detection power for a test using $p_0$ as critical p-value,
we need 
\begin{equation} \label{Gp0}
G(p_0) \geq p_0 + C \sqrt{p_0/N} \mbox{ for some } C>0. 
\end{equation}
Setting $C = \tfrac{2 \lambda_2}{\sqrt{\log N}}$ to the right-hand side of (\ref{Gp0}) gives us a function close to $F$.

As $\sqrt{\log N}$ varies slowly with $N$,  
we can view $F$ as lying at the boundary for which detection is possible at each small $p_0$.
This suggests that for any $G$ that is greater than $F$ at some $p_0$,
we are able to differentiate it from Uniform (0,1) by using the likelihood test against $F$.

When $p$ is of order smaller than $N^{-1}$,
$\tfrac{\log N}{N} f_1(p)$ dominates $\tfrac{1}{\sqrt{N \log N}} f_2(p)$
and the selection of $\lambda_1 > 0$ is advantageous.
This is relevant in the extension of sparse likelihood scores to detect change-points on long sequences where
large number of likelihood comparisons is involved.

\begin{figure}
\begin{center}
\includegraphics[width=2.5in]{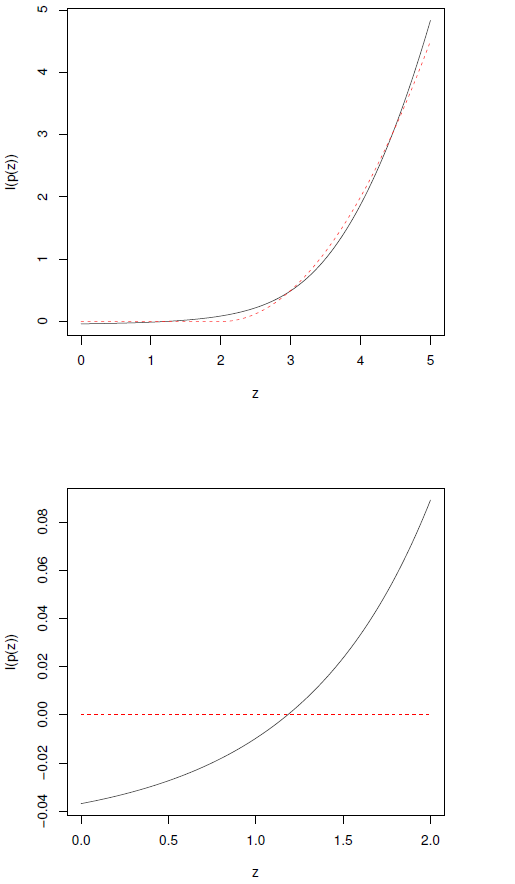}
\caption{ Graphs of $\ell_N(p(z))$ (black, ---) and $(z-2)^2_+/2$ (red,$--$),
with $p(z)=2 \Phi(-z)$,
for $0 \leq z \leq 5$ (top) and $0 \leq z \leq 2$ (bottom).
The parameters of $\ell_N$ are $N=500$,
$\lambda_1=1$ and $\lambda_2 = 1.84 \Big(\doteq \sqrt{\tfrac{\log T}{\log \log T}}$ for $T=500$\Big).}
\end{center}
\end{figure}

The sparse likelihood score can be viewed as a form of soft thresholding.
To visualize this we compare in Figure 1
the plot of $\ell_N(p(z))$ for $p(z)=2 \Phi(-|z|)$,
$N=500$, 
$\lambda_1=1$ and $\lambda_2=1.84$,
against that of $(z-2)^2_+/2$. 
For $0 \leq z \leq 5$,
the two functions are close to each other however within $0 \leq z \leq 2$,
$\ell_N(p(z))$ is not constant but has a gentle upward curve.
The sparsity likelihood score is negative for $z \leq 1.18$ and $\ell_N(p(Z))$ for $Z$ standard normal has a mean of $-0.004$. 
This negative mean helps in controlling the sum of scores when $N$ is large and $p^n \stackrel{\rm i.i.d}{\sim} {\rm Uniform}(0,1)$.

\subsection{Optimal detection}

We show here that the sparse likelihood score is optimal in the detection of change-points for a broad range of sparsity. 
Let $E_0$ and $P_0$ denote expectation and probability respectively with respect to $p^n \stackrel{\rm i.i.d}{\sim} {\rm Uniform}(0,1)$.
Since
\begin{eqnarray*}
& & E_0 \exp(\ell_N(\bp)) \cr
& = & \prod_{n=1}^N E_0[1+\tfrac{\lambda_1 \log N}{N} f_1(p^n)+\tfrac{\lambda_2}{\sqrt{N \log N}} f_2(p^n)]  = 1,
\end{eqnarray*}
it follows from Markov's inequality that
\begin{equation} \label{Mark} 
P_0(\ell_N(\bp) \geq c_N) \leq e^{-c_N}.
\end{equation}
This exponential bound makes the sparsity likelihood score easy to work with when there are large number of likelihood comparisons,
as critical values satisfying a required level of Type I error control can have a simple expression not depending on $N$.
We show in Theorem \ref{thm1} that by selecting
\begin{equation} \label{cN}
c_N \rightarrow \infty \mbox{ with } c_N =o(N^{\delta}) \rightarrow 0 \mbox{ for all } \delta>0,
\end{equation}
the Type I and II error probabilities both go to zero at the detection boundary.

\medskip
\begin{thm} \label{thm1}
Assume (\ref{cN}).
Consider the test of $H_0$: $Z^n \stackrel{\rm i.i.d.}{\sim} {\rm N}(0,1)$ versus $H_1$: $Z^n \stackrel{\rm i.i.d.}{\sim} (1-\epsilon) {\rm N}(0,1)+\epsilon {\rm N}(\mu_N,1)$,
for $1 \leq n \leq N$,
with $\epsilon = N^{-\beta}$ for some $\tfrac{1}{2} < \beta < 1$.
Let $p^n = \Phi(-Z^n)$.
If $\mu_N = \sqrt{2 \nu \log N}$ for $\nu > \rho_Z(\beta)$,
then 
$$P_0(\ell_N(\bp) \geq c_N)+P_{\mu_N}(\ell_N(\bp) < c_N) \rightarrow 0.
$$
\end{thm}

\section{Change-point detection}

Let $X^n_t$ denote the $t$th observation of the $n$th sequence for $1 \leq t \leq T$ and $1 \leq n \leq N$.
Consider first the model
\begin{equation} \label{normal}
X^n_t \sim_{\rm indep.} {\rm N}(\mu_t^n,1).
\end{equation}
We are interested in the detection and estimation of 
$$\btau := \{ t: \mu^n_t \neq \mu^n_{t+1} \mbox{ for some } n \}.
$$

For $s<t$,
let $\bar X_{st}^n = (t-s)^{-1} \sum_{u=s+1}^t X_u^n$.
To check for a change of mean on the $n$th sequence at location $t$,
select $s < t < u$ and let p-value
$$p_{stu}^n = 2 \Phi(-|Z_{stu}^n|), \mbox{ where } Z_{stu}^n = \tfrac{\bar X_{tu}^n - \bar X_{st}^n}{\sqrt{(u-t)^{-1}+(t-s)^{-1}}}.
$$
In the sparse likelihood algorithm we combine these p-values using $\ell_N(\bp_{stu})$,
where $\bp_{stu} = (p_{stu}^1, \ldots, p_{stu}^N)$.
When the data follow some other distributions, the corresponding likelihood ratio statistic and p-value can be computed accordingly.

Sparse likelihood scores detects well when only a small fraction of the sequences undergo change of mean.
For $T$ large
computing the sparse likelihood score for all $(s,t,u)$ is expensive. 
Instead we combine the approximating set idea of Walther (2010) to first space out the $(s,t,u)$ that are evaluated,
and to apply the CUSUM-type scores used in WBS to estimate the change-point location accurately only
when the first step indicates a change-point.  

In addition to computational savings,
through this two-step approach we are able to incorporate multi-scale penalization terms 
similar to those used in D\"{u}mbgen and Spokoiny (2001) and the SMUCE algorithm of Frick, Munk and Sieling (2014), 
to ensure optimality not only at all levels of sparse change-points,
but also at all orders of change magnitudes.

Let $1 \leq h_1 < h_2 < \cdots$ and $1 \leq d_1 < d_2 < \cdots$ be integer-valued sequences with $h_i \geq d_i$ for all $i$.
Let $K_i(g) = \lfloor \tfrac{g-1}{d_i} \rfloor$.
Define
\begin{eqnarray*}
\cA_i(g) & = & \{ (s(ik),t(ik),u(ik)): 1 \leq k \leq K_i(g) \}, \cr
s(ik) & = & \max(0,kd_i-h_i), \cr
t(ik) & = & kd_i, \cr
u(ik) & = & \min(kd_i+h_i,g).
\end{eqnarray*}
The elements of $\cA_i(g)$ are the indices where sparse likelihood scores for windows of length $h_i$ are computed.
Initially we have the full dataset $\bX_{1:T} = (X^n_t: 1 \leq t \leq T, 1 \leq n \leq N)$
and after one or more change-points have been estimated, 
it is split into sub-datasets $\bX_{b:e} = (X^n_t: b \leq t \leq e, 1 \leq n \leq N)$,
with length $g=e-b+1$.
We check for change-points in $\bX_{b:e}$ using windows specified by ${\cal A}_i(g)$.

Let the penalized sparse likelihood scores
\begin{equation} \label{plp}
\ell^{\rm{pen}}_N(\bp_{stu}) = \ell_N(\bp_{stu}) - \log(\tfrac{T}{4}(\tfrac{1}{t-s}+\tfrac{1}{u-t})).
\end{equation}
For $g \geq 1$,
let $i_g = \max \{ i: h_i+d_i \leq g \}$.
The detection of change-points within $\bX_{b:e}$,
with window lengths at least $h_{i_0}$, 
is as follows.

\begin{algorithm}
\caption{SL-estimate}
\begin{algorithmic}
\STATE 
\STATE \textsc{input}($c,i_0,b,e$)
\STATE \hspace{0.25cm} $\bX \leftarrow \bX_{b:e}$ 
\STATE \hspace{0.25cm} $g \leftarrow e-b+1$ 
\STATE \hspace{0.25cm} \textsc{for} $i=i_0, \ldots,i_g$
\STATE \hspace{0.5cm} \textsc{if} $\max_{1 \leq k \leq K_i(g)} \ell^{\rm{pen}}_N(\bp_{s(ik),t(ik),u(ik)}) \geq c$ \textsc{then}
\STATE \hspace{0.75cm} $j \leftarrow \mbox{argmax}_{k:1 \leq k \leq K_i(g)} \ell^{\rm{pen}}_N(\bp_{s(ik),t(ik),u(ik)})$
\STATE \hspace{0.75cm} $\wht \tau \leftarrow [\mbox{argmax}_{t:s(ij) < t < u(ij)} \ell^{\rm{pen}}_N(\bp_{s(ij),t,u(ij)})]+b-1$
\STATE \hspace{0.75cm} \textsc{output} $(\wht \tau,i)$
\STATE \hspace{0.75cm} \textsc{stop}
\STATE \hspace{0.5cm} \textsc{end if}
\STATE \hspace{0.25cm} \textsc{end for}
\STATE \hspace{0.25cm} \textsc{output} (0,0)
\end{algorithmic}
\end{algorithm}

There are two steps in SL-estimate in the estimation of a change-point,
when the largest penalized score exceeds the critical value $c$.
The first is the identification of an interval $(s(ij),u(ij))$,
associated with the largest penalized score,
within which a change-point lies.
The second is the estimation of the change-point within this interval.
In the approximating set $\cA_i(g)$,
neighboring windows are located $d_i$ apart,
hence we are unable to estimate the change-points accurately in the first step.
Accurate estimation is carried out, 
with more intensive computations within $(s(ij),u(ij))$,
in the second step.
Since the second step is performed only after an interval has been identified as containing a change-point,
performing this two-step procedure saves computations in regions where scores are generally small and the likelihood of change-points is low. 

After a change-point has been identified,
we split the dataset into two and execute the same algorithm on each split dataset.
To avoid repetitive computations,
we start from the window length $h_{i_0}$ used in the evaluation of the change-point spltting the dataset,
instead of starting from the smallest window length $h_1$,
on the split datasets.
The use of a set of representative set of window-lengths $h_i$ for computational savings in change-point detection have
been proposed in Willsky and Jones \cite{WJ76}.
The recursive segmentation algorithm for the computation of the estimated change-point set $\wht{\btau}$ is given below, 
with initialization at $(c,1,1,T,\emptyset)$.

\begin{figure}
\begin{center}
\includegraphics[width=2.5in]{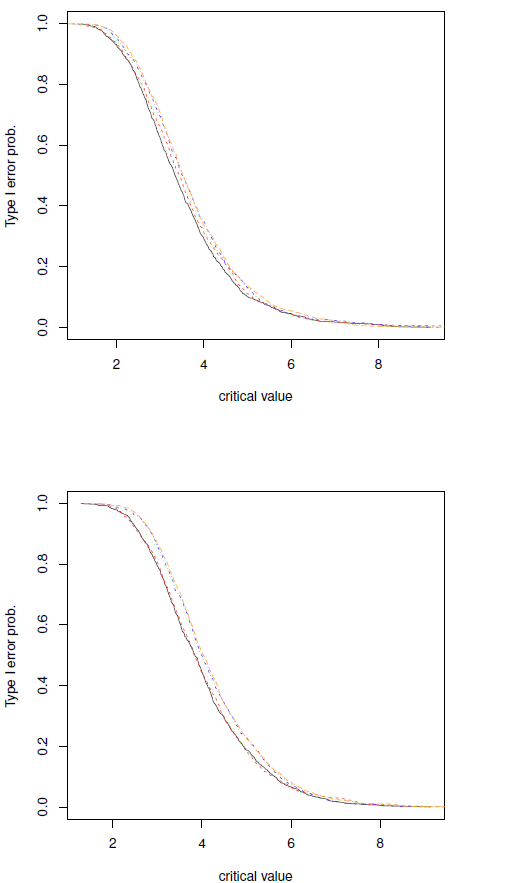}
\caption{ Graphs of Type I error probability against critical value for the sparse likelihood detection algorithm,
for independent unit variance normal observations. 
We consider parameters $d_i$, $h_i$, $\lambda_1$ and $\lambda_2$ as applied in the numerical studies in Section V,
with $T=2000$ (top), 
$T=$20,000 (bottom).
and $N=50$ (black), $N=100$ (red), $N=200$ (green), $N=500$ (blue), $N=1000$ (orange).} 
\end{center}
\end{figure}

\begin{algorithm}[H]
\caption{SL-detect}
\begin{algorithmic}
\STATE
\STATE \textsc{input}($c,i_0,b,e,\wht{\btau}$)
\STATE \hspace{0.5cm} $(\wht \tau,i) \leftarrow$ SL-estimate($c,i_0,b,e$) 
\STATE \hspace{0.5cm} \textsc{if} $\wht \tau>0$ \textsc{then}
\STATE \hspace{0.75cm} $\wht{\btau} \leftarrow \wht{\btau} \cup \{ \wht \tau \}$ 
\STATE \hspace{0.75cm} $\wht{\btau} \leftarrow$ SL-detect($c,i,b,\wht \tau,\wht{\btau})$ 
\STATE \hspace{0.75cm} $\wht{\btau} \leftarrow$ SL-detect($c,i,\wht \tau,e,\wht{\btau})$
\STATE \hspace{0.5cm} \textsc{end if}
\STATE \hspace{0.5cm} \textsc{output} $\wht{\btau}$ 
\end{algorithmic}
\end{algorithm}

In Figure 2 we show that the critical values of the sparse likelihood algorithm, 
for a specified Type I error probability,
is stable over $N$.
Contributing factors include $\ell_N(p)$ having a mean that is close to zero and $\ell_N(\bp)$ having exponential tail probabilities not depending on $N$,
see (\ref{Mark}),
when $p$ and $p^n$ are uniformly distributed.

\section{Optimal detection}

Let $\bmu = (\mu_t^n: 1 \leq t \leq T, 1 \leq n \leq N)$ and let $J= (\# \btau)$ be the number of change-points.
We show that the sparse likelihood algorithm is optimal for normal observations in Section IV.A
and for Poisson observations in Section IV.B.
Consider $T \rightarrow \infty$ as $N \rightarrow \infty$ such that
\begin{equation} \label{zeta}
\log T \sim N^{\zeta} \mbox{ for some } 0 < \zeta < 1.
\end{equation}

In Theorems \ref{thm2} and \ref{thm4} 
we specify the detection boundary for asymptotically zero Type I and II error probabilities. 
Analogous detection boundaries for a single sequence is given in Arias-Castro, Donoho and Huo \cite{ADH05, ADH06}. 

In Theorems \ref{thm3} and \ref{thm5} we show that Type I and II error probabilities of the sparse likelihood algorithm go to zero at the detection boundary.

Recall that $i_T = \max \{ i: h_i+d_i \leq T \}$.
Consider the sparse likelihood algorithm with $d_i$ and $h_i$ satisfying
\begin{eqnarray} \label{hi}
\tfrac{h_{i+1}}{h_i} & \rightarrow & 1 \mbox{ and } d_i = o(h_i) \mbox{ as } i \rightarrow \infty, \\ \label{iT}
\log \Big( \sum_{i=1}^{i_T} \tfrac{h_i}{d_i} \Big) & = & o(\log T) \mbox{ as } T \rightarrow \infty,
\end{eqnarray}
and critical values $c_T$ satisfying
\begin{equation} \label{cTT}
c_T = o(\log T) \mbox{ and } c_T - \log \Big( \sum_{i=1}^{i_T} \tfrac{h_i}{d_i} \Big) \rightarrow \infty \mbox{ as } T \rightarrow \infty.
\end{equation}
For the sparse likelihood algorithm select parameters $\lambda_1 > 0$ and
\begin{equation} \label{lam12}
\lambda_2 = \sqrt{\tfrac{\log T}{\log \log T}}.
\end{equation}

We satisfy (\ref{hi}) when $h_i \sim \exp(\tfrac{i}{\log i})$ and $d_i \sim \tfrac{h_i}{i}$ as $i~\rightarrow~\infty$.
Moreover (\ref{iT}) holds because
$$\log \Big( \sum_{i=1}^{i_T} \tfrac{h_i}{d_i} \Big) \sim 2 \log i_T \sim 2 \log \log T.
$$

Condition (\ref{hi}) ensures that the set of $(h_i,d_i)$ is sufficiently dense to detect change-points optimally.
Condition (\ref{iT}) is required for (\ref{cTT}) to hold.
The first half of condition (\ref{cTT}) ensures Type II error probability goes to 0.
The second half ensures Type I error probability goes to 0.



\subsection{Normal model}

Let 
$$m_{j \Delta} = \# \{ n: |\mu_{\tau_j+1}^n-\mu_{\tau_j}^n| \geq \Delta \}
$$
be the number of sequences with change of mean of at least $\Delta$ at the $j$th change-point.
Let
\begin{eqnarray*}
\Omega_0 & = & \{ \bmu: J=0 \}, \cr
\Omega_1(\Delta,V,h) & = & \{ \bmu: \mbox{ there exists } j \mbox{ such that} \cr
& & \qquad \min(\tau_j-\tau_{j-1},\tau_{j+1}-\tau_j) \geq h \cr
& & \qquad \mbox{ and } m_{j \Delta} \geq V \}, 
\end{eqnarray*}
with the convention $\tau_0=0$ and $\tau_{J+1}=T$.
We consider here the test of $H_0$: $\bmu \in \Omega_0$ versus $H_1$: $\bmu \in \Omega_1(\Delta,h,V)$.
Define
\begin{equation} \label{rhoZb}
\rho_Z(\beta,\zeta) = \left\{ \begin{array}{l} \beta-\tfrac{1-\zeta}{2} \cr
\qquad \mbox{ if } \tfrac{1-\zeta}{2} < \beta \leq 
\tfrac{3(1-\zeta)}{4}, \cr
(\sqrt{1-\zeta}-\sqrt{1-\zeta-\beta})^2 \cr
\qquad \mbox{ if } \tfrac{3(1-\zeta)}{4} < \beta < 1-\zeta. \end{array} \right.
\end{equation}
These constants are extensions of $\rho_Z(\beta)$ in (\ref{rhoZ}) to capture the effect of multiple testing in change-point detection.

\medskip
\begin{thm} \label{thm2}
Assume {\rm (\ref{zeta})} and let $0 < \epsilon < 1$.
For normal observations,
no algorithm is able to achieve,
as $N  \rightarrow \infty$,
\begin{equation} \label{zomega}
\sup_{\bmu \in \Omega_0} P_{\bmu}(\mbox{Type I error})+\sup_{\bmu \in \Omega_1(\Delta,V,h)} P_{\bmu}(\mbox{Type II error}) \rightarrow 0,
\end{equation}
under any of the following conditions.

\medskip
{\rm (a)} Consider $\Delta>0$ fixed.

\smallskip
\qquad {i.} When $V=o(\tfrac{\log T}{\log N})$ and $h = 4 (1-\epsilon)(\tfrac{\log T}{\Delta^2 V})$.

\smallskip
\qquad {\rm ii.} When $V \sim N^{1-\beta}$ for some $\tfrac{1-\zeta}{2} < \beta < 1-\zeta$ and 
$h = 4(1-\epsilon) \rho_Z(\beta,\zeta) (\tfrac{\log N}{\Delta^2})$.

\medskip
{\rm (b)} Consider $\Delta = T^{-\eta}$ for some $0 < \eta < \tfrac{1}{2}$.

\smallskip
\qquad {\rm i.} When $V=o(\tfrac{\log T}{\log N})$ and $h = 4(1-2 \eta)(1-\epsilon)(\tfrac{\log T}{\Delta^2 V})$. 

\smallskip
\qquad {\rm ii.} When $V \sim N^{1-\beta}$ for some $\tfrac{1-\zeta}{2} < \beta < 1-\zeta$ and 
$h = 4(1-\epsilon) \rho_Z(\beta,\zeta) (\tfrac{\log N}{\Delta^2})$.
\end{thm}

\medskip
\begin{thm} \label{thm3}
Assume {\rm (\ref{zeta})} and let $\epsilon>0$.
For normal observations the sparse likelihood algorithm, 
with parameters satisfying {\rm (\ref{hi})--(\ref{lam12})}, 
achieves {\rm (\ref{zomega})} under any of the following conditions.

\medskip
{\rm (a)} Consider $\Delta>0$ fixed.

\smallskip
\qquad {\rm i.} When $V=o(\tfrac{\log T}{\log N})$ and $h = 4 (1+\epsilon)(\tfrac{\log T}{\Delta^2 V})$. 

\smallskip
\qquad {\rm ii.} When $V \sim N^{1-\beta}$ for some $\tfrac{1-\zeta}{2} < \beta < 1-\zeta$ and 
$h = 4(1+\epsilon) \rho_Z(\beta,\zeta)(\tfrac{\log N}{\Delta^2})$.

\medskip
{\rm (b)} Consider $\Delta = T^{-\eta}$ for some $0 < \eta < \tfrac{1}{2}$.

\smallskip
\qquad {\rm i.} When $V=o(\tfrac{N^{\zeta}}{\log N})$ and $h = 4(1-2 \eta)(1+\epsilon)(\tfrac{\log T}{\Delta^2 V})$.

\smallskip
\qquad {\rm ii.} When $V \sim N^{1-\beta}$ for some $\tfrac{1-\zeta}{2} < \beta < 1-\zeta$ and 
$h = 4(1+\epsilon) \rho_Z(\beta,\zeta) (\tfrac{\log N}{\Delta^2})$.
\end{thm}


\subsection{Poisson model}

We show here the optimality of the sparse likelihood detection algorithm for 
\begin{equation} \label{XP}
X^n_t \sim_{\rm indep.} \mbox{Poisson}(\mu^n_t).
\end{equation}
For optimal detection on a single Poisson sequence,
see Rivera and Walther \cite{RW13}.   

Let $Y_{st}^n = \sum_{v=s+1}^t X_v^n$.
Consider $s<t<u$.
Under the null hypothesis of no change-points in the interval $(s,u)$,
conditioned on $Y_{su}^n=y^n_{su}$,
$Y^n_{st}$ is binomial distributed with $y^n_{su}$ trials and success probability $\frac{t-s}{u-s}$.
Let $p^n_{stu}$ be the two-sided p-value of this conditional binomial test,
with continuity adjustments so that it is distributed as Uniform(0,1) under the null hypothesis.
More specifically when $Y_{st}^n = y_{st}^n$ and $Y_{su}^n = y_{su}^n$ simulate
\begin{equation} \label{psi}
\psi_{stu}^n \sim \mbox{Uniform}(P(Y < y_{st}^n),P(Y \leq y_{st}^n)),
\end{equation}
where $P$ is probability with respect to $Y \sim \mbox{Binomial}(y_{su}^n,\tfrac{t-s}{u-s})$, 
and define $p_{stu}^n = 2 \min(\psi_{stu}^n,1-\psi_{stu}^n)$. 

Let
$$m_{j \Delta} = \# \{ n: |\log (\mu_{\tau_j+1}^n/\mu_{\tau_j}^n)| \geq \Delta \},
$$
and for a given $\mu_0 > 0$,
let 
\begin{eqnarray*}
\Lambda & = & \{ \bmu: \mu_n^t \geq \mu_0 \mbox{ for all } n \mbox{ and } t \}, \cr
\Lambda_0 & = & \{ \bmu \in \Lambda: J = 0 \}, \cr
\Lambda_1(\Delta,V,h) & = & \{ \bmu \in \Lambda: \mbox{ there exists } j \mbox{ such that } \cr
& & \qquad \min(\tau_{j+1}-\tau_j,\tau_j-\tau_{j-1}) \geq h \cr
& & \qquad \mbox{ and } m_{j \Delta} \geq V \}.  
\end{eqnarray*}
We consider here the test of $H_0$: $\bmu \in \Lambda_0$ vs $H_1$: $\bmu \in \Lambda_1(\Delta,V,h)$.

For a given $r>1$, 
let 
\begin{equation} \label{IY}
I_r = r \log(\tfrac{2r}{r+1})+\log(\tfrac{2}{r+1}).
\end{equation}
Let $g_r(\omega) = (\tfrac{1+r^{\omega}}{2})^{\frac{1}{\omega}}$ and let 
\begin{equation} \label{rho}
\rho_r(\beta,\zeta) = \max_{\frac{1-\zeta}{\beta} < \omega \leq 2} \Big( \tfrac{\beta-\omega^{-1}(1-\zeta)}{2 g_r(\omega)-1-r} \Big)
\mbox{ for } \tfrac{1-\zeta}{2} < \beta < 1-\zeta.
\end{equation}

In Theorem \ref{thm4} we show that (\ref{rho}) is the asymptotic constant in the detection boundary of Poisson random variables.
In Theorem \ref{thm5} we show that the sparse likelihood algorithm achieves detection at this boundary for a broad
range of sparsity. 

\medskip
\begin{thm} \label{thm4}
Assume {\rm (\ref{zeta})}.
Let $r=e^{\Delta}$ for some $\Delta>0$ and $0 < \epsilon < 1$.
For Poisson observations no algorithm is able to achieve,
as $N \rightarrow \infty$,
\begin{equation} \label{zlambda}
\sup_{\bmu \in \Lambda_0} P_{\bmu}(\mbox{Type I error})+\sup_{\bmu \in \Lambda_1(\Delta,V,h)} P_{\bmu}(\mbox{Type II error}) \rightarrow 0
\end{equation}
under either of the following conditions.

\medskip
{\rm (a)} When $V=o(\tfrac{\log T}{\log N})$ and $h = (1-\epsilon) I_r^{-1} (\tfrac{\log T}{\mu_0 V})$. 

\smallskip
{\rm (b)} When $V \sim N^{1-\beta}$ for some $\tfrac{1-\zeta}{2} < \beta < 1-\zeta$ and 
$h = (1-\epsilon) \rho_r(\beta,\zeta) (\tfrac{\log N}{\mu_0})$. 
\end{thm}

\bigskip
\begin{thm} \label{thm5}
Assume {\rm (\ref{zeta})}.
Let $\epsilon > 0$, $\Delta>0$ and $1 < r < e^{\Delta}$.
For Poisson observations the sparse likelihood algorithm, 
with parameters satisfying {\rm (\ref{hi})--(\ref{lam12})},
achieves {\rm (\ref{zlambda})} under either of the following conditions.

\medskip
{\rm (a)} When $V=o(\tfrac{\log T}{\log N})$ and $h = (1+\epsilon) I_r^{-1} (\tfrac{\log T}{\mu_0 V})$.

\smallskip
{\rm (b)} When $V \sim N^{1-\beta}$ for some $\tfrac{1-\zeta}{2} < \beta < 1-\zeta$ and 
$h = (1+\epsilon) \rho_r(\beta,\zeta) (\tfrac{\log N}{\mu_0})$. 
\end{thm}



\section{Simulation studies}

We follow here the simulation set-up in Sections 5.1 and 5.3 of Wang and Samworth \cite{WS18}.
Assume that the random variables are normal with variances that are unknown but equal within sequence.
These variances are estimated using median absolute differences of adjacent observations and after normalization,
the random variables are treated like unit variance normal. 

In the first study there is exactly one change-point $\tau_1$.
Consider $\mu_t^n=0$ for $t \leq \tau_1$ and all $n$.
For $t>\tau_1$, 
let
$$\mu^t_n = \left\{ \begin{array}{ll} 0.8 \Big/ \sqrt{n \textstyle\sum_{m=1}^V m^{-1}} & \mbox{ if } n \leq V, \cr
0 & \mbox{ if } n>V. \end{array} \right.
$$
The objective is to estimate $\tau_1$ assuming we know there is exactly one change-point.
We estimate $\tau_1$ here by
$$\wht \tau_1 = {\rm arg} \max_{0 < t < T} \ell^{\rm{pen}}_N(\bp_{0tT}),
$$
where $\ell^{\rm pen}_N$ is the penalized sparse score with $\lambda_1=1$ and $\lambda_2 = \sqrt{\tfrac{\log T}{\log \log T}}$.

We simulate the probabilities that $|\wht \tau_1-\tau_1| \leq k$ for $k=3$ and 10,
and compare against the INSPECT algorithm and the scan algorithm of Enikeeva and Harchaoui \cite{EH19}.
These two algorithms have the best numerical performances in Wang and Samworth \cite{WS18}.
The comparisons in Table~I show that the sparse likelihood algorithm performs well.  

\begin{table}[!t]
\caption{The fraction of simulation runs (out of 1000) for which $\wht \tau_1$ is within distance $k$ from $\tau_1$ for $k=3$ and 10.
The same datasets are used to compare sparse likelihood (SL), INSPECT and the scan test,
with $\tau_1=200$ for $T=500$ and $\tau_1=800$ for $T=2000$.}
\centering
\begin{tabular}{rr|rr|rr|rr}
& $k$ & 3 & 10 & 3 & 10 & 3 & 10 \cr 
& $V$ & \multicolumn{2}{c|}{SL}  & \multicolumn{2}{c|}{INSPECT} & \multicolumn{2}{c}{scan} \cr \hline
$T=500$ & 3 &  0.511 & 0.801 & 0.478 & 0.785 & {\it 0.520} & {\it 0.804} \cr
$N=500$ & 5 & {\it 0.466} & {\it 0.740} & 0.427 & 0.718 & 0.463 & 0.722 \cr
& 10 & {\it 0.393} & {\it 0.645} & 0.370 & 0.637 & 0.362 & 0.599 \cr
& 22 & {\it 0.319} & {\it 0.553} & 0.282 & 0.547 & 0.256 & 0.465 \cr
& 50 & {\it 0.244} & {\it 0.462} & 0.211 & 0.453 & 0.197 & 0.378 \cr
& 500 & {\it 0.177} & {\it 0.339} & 0.148 & 0.335 & 0.112 & 0.240 \cr \hline
$T=500$ & 3 & {\it 0.481} & {\it 0.748} & 0.410 & 0.667 & 0.480 & 0.730 \cr
$N=2000$ & 5 & {\it 0.423} & {\it 0.673} & 0.344 & 0.584 & 0.394 & 0.633 \cr
& 10 & {\it 0.320} & {\it 0.546} & 0.246 & 0.480 & 0.261 & 0.456 \cr
& 20 & {\it 0.237} & {\it 0.431} & 0.198 & 0.403 & 0.188 & 0.332 \cr
& 45 & {\it 0.186} & {\it 0.344} & 0.136 & 0.311 & 0.130 & 0.242 \cr
& 200 & {\it 0.114} & 0.227 & 0.095 & {\it 0.235} & 0.074 & 0.153 \cr
& 2000 & 0.068 & 0.160 & {\it 0.078} & {\it 0.189} & 0.042 & 0.096 \cr \hline
$T=2000$ & 3 & {\it 0.603} & {\it 0.859} & 0.587 & 0.855 & 0.589 & 0.854 \cr
$N=500$ & 5 & {\it 0.604} & {\it 0.865} & 0.595 & 0.855 & 0.558 & 0.832 \cr
& 10 & 0.565 & 0.827 & {\it 0.569} & {\it 0.833} & 0.487 & 0.764 \cr
& 22 & {\it 0.522} & 0.789 & {\it 0.522} & {\it 0.795} & 0.438 & 0.714 \cr
& 50 & {\it 0.472} & {\it 0.748} & 0.468 & 0.745 & 0.384 & 0.652 \cr
& 500 & {\it 0.378} & {\it 0.643} & 0.336 & 0.609 & 0.273 & 0.524 \cr \hline
$T=2000$ & 3 & 0.607 & {\it 0.866} & {\it 0.608} & 0.861 & 0.599 & 0.858 \cr
$N=2000$ & 5 & {\it 0.594} & {\it 0.864} & 0.586 & 0.857 & 0.557 & 0.829 \cr
& 10 & 0.553 & {\it 0.847} & {\it 0.558} & {\it 0.847} & 0.476 & 0.780 \cr
& 20 & 0.494 & {\it 0.807} & {\it 0.498} & 0.789 & 0.435 & 0.726 \cr
& 45 & 0.447 & {\it 0.747} & {\it 0.451} & 0.746 & 0.377 & 0.657 \cr
& 200 & {\it 0.362} & {\it 0.649} & 0.342 & 0.604 & 0.297 & 0.554 \cr
& 2000 & {\it 0.274} & {\it 0.532} & 0.241 & 0.471 & 0.225 & 0.457 \cr
\end{tabular}
\end{table}

In the second study there are three change-points within $N=200$ sequences of length $T=2000$,
at $\tau_1=500$, $\tau_2=1000$ and $\tau_3=1500$.
At each change-point exactly $40$ sequences undergo mean changes.
Six scenarios are considered,
corresponding to
$$\mu_{\tau_j+1}^{k(j-1)+n}-\mu_{\tau_j}^{k(j-1)+n} = r \Big/ \sqrt{n \textstyle\sum_{m=1}^{40} m^{-1}},
$$
$1 \leq j \leq 3$ and $1 \leq n \leq 40$, 
for $r=0.4, 0.6$ and $k=0, 20, 40$.
For $k=0$, 
the mean changes are within the same 40 sequences at all three change-points,
whereas for $k=40$ the mean changes at all three change-points are on distinct sequences. 
For $k=20$,
there is partial overlap of the sequences having mean changes at adjacent change-points.
The number of estimated change-points over 100 simulated datasets on each sequence is recorded,
as well as the adjusted Rand index (ARI),
see Rand \cite{Ran71} and Hubert and Arabie \cite{HA85}, 
to measure the quality of the change-point estimation.

In the application of the sparse likelihood algorithm,
we select $h_1=1$ and $h_{i+1}=\lceil 1.1 h_i \rceil$ for $i \geq 1$,
and $d_i = \lfloor h_i/i \rfloor$,
for a total of $i_T=61$ window lengths.
We select critical value $c_T=5$ and parameters $\lambda_1=1$, 
$\lambda_2 = \sqrt{\tfrac{\log T}{\log \log T}} \doteq 1.94$.

Wang and Samworth \cite{WS18} showed that INSPECT achieves average ARI of 0.90 when $r=0.6$ and either 0.73 (for $k=20$) or 0.74 
(for $k=0$ and 40) when $r=0.4$,
comparable to sparse likelihood,
see Table II.

In addition to INSPECT,
Wang and Samworth \cite{WS18} considered DC, SBS and scan, 
as well as the CUSUM aggregration algorithms of Jirak \cite{Jir15} and Horv\'{a}th and Hu\v{s}kov\'{a} \cite{HH14}, 
with average ARI in the range 0.77--0.87 when $r=0.6$ and 0.68--0.72 when $r=0.4$.

\begin{table}[!t]
\caption{Number of change-points estimated by the sparse likelihood algorithm and the average ARI over 100 simulated datasets.}
\centering
\begin{tabular}{rr|rrrr|r}
$r$ & $k$ & \multicolumn{4}{c|}{$\#$ change-points} & ARI \cr
& & 2 & 3 & 4 & 5 \cr \hline
0.6 & 0 & 11 & 80 & 8 & 1 & 0.91 \cr
0.4 & 0 & 61 & 35 & 4 & 0 & 0.74 \cr
0.6 & 20 & 12 & 80 & 8 & 0 & 0.91 \cr
0.4 & 20 & 66 & 31 & 2 & 1 & 0.74 \cr
0.6 & 40 & 10 & 78 & 12 & 0 & 0.91 \cr
0.4 & 40 & 68 & 26 & 6 & 0 & 0.75 \cr
\end{tabular}
\end{table}

\begin{appendix}
\section{Proof of Theorem \ref{thm1}}

Since $c_N \rightarrow \infty$,
by Markov's inequality $P_0(\ell_N(\bp) \geq c_N) \leq e^{-c_N} \rightarrow 0$.
The proof of $P_{\mu_N}(\ell_N(\bp) < c_N) \rightarrow 0$ applies Lemmas \ref{lem1} and \ref{lem2} below.
Lemma \ref{lem1} says that the sum of sparse likelihood scores under $q^n \sim {\rm Uniform}(0,1)$ is bounded below
by a value close to zero,
with large probability.
Lemma \ref{lem2} provides a lower bound to the increase in score when the p-value is divided by at least 2.
Their proofs are at the end of Appendix~A.

\medskip
\begin{lem} \label{lem1}
Let $\bq = (q^1, \ldots, q^N)$,
with $q^n \sim_{\rm i.i.d.}$ {\rm Uniform(0,1)}.
For fixed $\lambda_1 \geq 0$ and $\delta > 0$,
$$\sup_{\delta \leq \lambda_2 \leq \sqrt{N}} P ( \ell_N(\bq) \leq -\lambda_2^2 \sqrt{\log N}) \rightarrow 0.
$$
\end{lem}

\medskip
\begin{lem} \label{lem2}
For $\lambda_1 > 0$ fixed, 
$\delta \leq \lambda_2 \leq \sqrt{N}$ for some $\delta>0$ and $\xi_N = o(N^{-\eta})$ for some
$\eta>0$ such that $\xi_N \geq \tfrac{\lambda_2^2}{2N}$,
$$\inf_{\substack{(p,q): p \leq \xi_N, \\ q \geq \lambda_2^2/N, p \leq q/2}} [\ell_N(p)-\ell_N(q)]
\geq \tfrac{\lambda_2}{4 \sqrt{N \xi_N \log N}}
$$
for large $N$.
\end{lem}

\medskip
{\sc Proof of Theorem} \ref{thm1}. 
Let $\tfrac{1}{2} < \beta < 1$, $\lambda_1 \geq 0$ and $\lambda_2>0$ be fixed.
Let $\nu$ be such that
\begin{eqnarray*} 
(1-\sqrt{1-\beta})^2 < \nu < 1 & &
\mbox{ if } \tfrac{3}{4} \leq \beta < 1, \cr
\beta-\tfrac{1}{2} < \nu < 4(\beta-\tfrac{1}{2}) & & \mbox{ if } \tfrac{1}{2} < \beta < \tfrac{3}{4}, 
\end{eqnarray*} 
and let
\begin{eqnarray*}
\mu_N & = & \sqrt{2 \nu \log N}, \cr
Q^n & \sim & \mbox{Bernoulli}(N^{-\beta}), \cr
Z^n|Q^n & \sim & {\rm N}(\mu_N Q^n,1), \cr
p^n & = & \Phi(-Z^n), \cr
q^n & = & \Phi(-Z^n+\mu_N Q^n).
\end{eqnarray*}
The additional assumptions of $\nu<1$ for $\tfrac{3}{4} \leq \beta < 1$ and $\nu < 4(\beta-\tfrac{1}{2})$ for $\tfrac{1}{2} < 
\beta < \tfrac{3}{4}$ is not restrictive because $\ell_N(\bp)$ increases stochastically with $\mu_N$.

\medskip
{\bf Case 1}: $\tfrac{3}{4} \leq \beta < 1$.
Let
$$\Gamma  =\{ n: Q^n=1, Z^n \geq \sqrt{2 \log N}, q^n \geq \tfrac{\lambda^2_2}{N} \}.
$$
For $N$ large,
$p^n \leq \tfrac{q^n}{2}$ for $n \in \Gamma$.
Moreover $\ell_N(p^n) \geq \ell_N(q^n)$ for all $n$.
Hence by Lemma \ref{lem1} and Lemma \ref{lem2} with $\xi_N = N^{-1}$,
with probability tending to 1,
\begin{eqnarray*}
\ell_N(\bp) & \geq & \ell_N(\bq) + \sum_{n \in \Gamma} [\ell_N(p^n)-\ell_N(q^n)] \cr
& \geq & -\lambda_2^2 \sqrt{\log N} + (\# \Gamma) \tfrac{\lambda_2}{4 \sqrt{\log N}}.
\end{eqnarray*}

Since $\# \Gamma$ is binomial with mean
\begin{eqnarray} \label{rN}
& & E_{\mu_N}(\# \Gamma) \\ \nonumber
& = & N^{1-\beta} [\Phi(-\sqrt{2 \log N}+\sqrt{2 \nu \log N}) -\tfrac{\lambda_2^2}{N}] \\ \nonumber
& \gtrsim & \tfrac{N^{1-\beta-(1-\sqrt{\nu})^2}}{\sqrt{\log N}},
\end{eqnarray}
with  $1-\beta-(1-\sqrt{\nu})^2 > 0$ for $(1-\sqrt{1-\beta})^2 < \nu < 1$, 
and since $c_N$ is subpolynomial in $N$,
we conclude $P_{\mu_N}(\ell_N(\bp) \geq c_N) \rightarrow 1$.

\medskip
{\bf Case 2}: $\tfrac{1}{2} < \beta < \tfrac{3}{4}$.
Let
$$\Gamma = \{ n: Q^n=1, Z^n \geq 2 \sqrt{(2 \beta-1) \log N}, q^n \geq \tfrac{\lambda_2^2}{N} \}. 
$$
For $N$ large,
$p^n \leq \tfrac{q^n}{2}$ for $n \in \Gamma$.
Hence by Lemma \ref{lem1} and Lemma \ref{lem2} with $\xi_N = N^{2-4 \beta}$,
with probability tending to 1,
\begin{eqnarray*}
\ell_N(\bp) & \geq & \ell_N(\bq) + \sum_{n \in \Gamma} [\ell_N(p^n)-\ell_N(q^n)] \cr
& \geq & -\lambda_2^2 \sqrt{\log N} + (\# \Gamma) \tfrac{\lambda_2}{4 N^{\frac{3}{2}-2 \beta} \sqrt{\log N}}.
\end{eqnarray*}

Since $\# \Gamma$ is binomial with mean
\begin{eqnarray} \label{rN}
& & E_{\mu_N}(\# \Gamma) \\ \nonumber
& = & N^{1-\beta} [\Phi(-2 \sqrt{(2 \beta-1) \log N}+\sqrt{2 \nu \log N}) -\tfrac{\lambda_2^2}{N}] \\ \nonumber
& \gtrsim & \tfrac{N^{1-\beta-(\sqrt{4 \beta-2}-\sqrt{\nu})^2}}{\sqrt{\log N}},
\end{eqnarray}
and
$$1-\beta-(\sqrt{4 \beta-2}-\sqrt{\nu})^2 > \tfrac{3}{2}-\beta \mbox{ for } \beta-\tfrac{1}{2} < \nu < 4(\beta-\tfrac{1}{2}),
$$ 
we conclude
$P_{\mu_N}(\ell_N(\bp) \geq c_N) \rightarrow 1$. 
$\wbox$

\medskip
{\sc Proof of Lemma} \ref{lem1}.
Let
$$x_N(p) = \tfrac{\lambda_1 \log N}{N} f_1(p)+\tfrac{\lambda_2}{\sqrt{N \log N}} f_2(p),
$$
where $f_1(p)=\tfrac{1}{p(2-\log p)^2}-\frac{1}{2}$, 
$f_2(p)=\tfrac{1}{\sqrt{p}}-2$,
$\lambda_1 \geq 0$ and $\delta \leq \lambda_2 \leq N^{\frac{1}{2}}$ for some $\delta>0$.
Let $r_N = \tfrac{1}{N \log N}$.
Since $x_N(r_N) \geq 0$ and $x_N(1) \geq -\tfrac{1}{2}$ for $N$ large and 
$\log(1+x) \geq x-x^2$ for $x \geq -\frac{1}{2}$, 
\begin{equation} \label{ellq}
\ell_N(\bq) = \sum_{n=1}^N \log(1+x_N(q^n)) \geq \sum_{n=1}^N h_N(q^n) - \sum_{n=1}^N h_N^2(q^n),
\end{equation}
where $h_N(q) = x_N(q) {\bf 1}_{\{ q \geq r_N \}}$.

By Chebyshev's inequality and the bounds in (\ref{in1})--(\ref{in3}) below.
\begin{eqnarray*}
& & P(\ell_N(\bq) \leq -\lambda_2^2 \sqrt{\log N}) \cr
& \leq & P \Big( \sum_{n=1}^N h_N(q^n) \leq -\tfrac{\lambda_2^2 \sqrt{\log N}}{2} \Big) \cr
& & \quad + P \Big( \sum_{n=1}^N h_N^2(q^n) \geq \tfrac{\lambda_2^2 \sqrt{\log N}}{2} \Big) \cr
& \leq & \tfrac{N \var(h_N(q^n))}{(N E h_N(q^n)+\frac{\lambda_2^2 \sqrt{\log N}}{2})^2}
+ \tfrac{N \var(h_N^2(q^n))}{(\frac{\lambda_2^2 \sqrt{\log N}}{2}-N E h_N^2(q^n))^2} \rightarrow 0.
\end{eqnarray*}

\medskip
Since $E x_N(q^n)=0$,
\begin{eqnarray} \label{in1}
E h_N(q^n) & = & - E[x_N(q^n) {\bf 1}_{\{ q^n < r_N \}}] \\ \nonumber
& = & -\tfrac{\lambda_1 \log N}{N}(\tfrac{1}{2-\log r_N}-\tfrac{r_N}{2}) \\ \nonumber
& & \quad -\tfrac{\lambda_2}{\sqrt{N \log N}}(2 \sqrt{r_N}-2 r_N) \\ \nonumber
& \geq & - \tfrac{\lambda_1}{N}-\tfrac{2 \lambda_2}{N \log N}.
\end{eqnarray}

\medskip
Let $s_N = \tfrac{(\log N)^2}{N}$.
\begin{eqnarray} \label{in2}
& & \var(h_N(q^n)) \\ \nonumber
& \leq & E h_N^2(q^n) \\ \nonumber
& \leq & \tfrac{2 \lambda_1^2 (\log N)^2}{N^2} \int_{r_N}^1 \tfrac{dq}{q^2 (2-\log q)^4} 
+ \tfrac{2 \lambda_2^2}{N \log N} \int_{r_N}^1 \tfrac{dq}{q} \\ \nonumber
& \leq & \tfrac{2 \lambda_1^2 (\log N)^2}{N^2} \Big( \int_{s_N}^1 \tfrac{dq}{q^2} + \tfrac{1}{(2-\log s_N)^4} \int_{r_N}^{s_N}
\tfrac{dq}{q^2}\Big) \\ \nonumber
& & \quad + \tfrac{2 \lambda_2^2 \log(\frac{1}{r_N})}{N \log N} \\ \nonumber
& \lesssim & \tfrac{\lambda_1^2+\lambda_2^2}{N}.
\end{eqnarray}

\medskip 
\begin{eqnarray} \label{in3}
& & \var(h_N^2(q^n)) \\ \nonumber
& \leq & E h_N^4(q^n) \\ \nonumber
& \leq & \tfrac{8 \lambda_1^4 (\log N)^4}{N^4} \Big( \int_{s_N}^1 \tfrac{dq}{q^4} + \tfrac{1}{(2-\log s_N)^8}
\int_{r_N}^{s_N} \tfrac{dq}{q^4} \Big) \\ \nonumber
& & \quad + \tfrac{8 \lambda_2^4}{(N \log N)^4} \int_{r_N}^1 \tfrac{dq}{q^2} \\ \nonumber
& \lesssim & \tfrac{\lambda_1^4+\lambda_2^4}{N}. 
\end{eqnarray}
$\wbox$

\medskip
{\sc Proof of Lemma} \ref{lem2}.
For $\tfrac{\lambda_2^2}{2N} \leq r \leq 2 \xi_N$,
$|\log r| \asymp \log N$ and therefore 
$$\tfrac{\frac{\lambda_1 \log N}{N} f_1(r)}{\frac{\lambda_2}{\sqrt{N \log N}} f_2(r)} \asymp  
\tfrac{1}{\lambda_2 \sqrt{Nr \log N}} \rightarrow 0.
$$
Moreover,
$$\tfrac{\lambda_2}{N \log N} f_2(r) \sim \tfrac{\lambda_2}{\sqrt{Nr \log N}} \rightarrow 0.
$$
Hence by $\log(1+x) \sim x$ as $x \rightarrow 0$,
\begin{equation} \label{lapprox}
\ell_N(r) \sim \tfrac{\lambda_2}{\sqrt{Nr \log N}}.
\end{equation} 

\medskip
{\bf Case 1}: $\frac{\lambda_2^2}{2N} \leq p \leq \xi_N$.
By (\ref{lapprox}) and $q \geq 2p$,
\begin{eqnarray*}
\ell_N(p)-\ell_N(q) & \geq & \ell_N(p)-\ell_N(2p) \cr
& \sim & (1-\tfrac{1}{\sqrt{2}}) \tfrac{\lambda_2}{\sqrt{Np \log N}} \cr
& > & \tfrac{\lambda_2}{4 \sqrt{N \xi_N \log N}}.
\end{eqnarray*}

\medskip
{\bf Case 2}: $p < \tfrac{\lambda_2^2}{2N}$.
By (\ref{lapprox}),
$q \geq \tfrac{\lambda_2^2}{N}$ and $\xi_N \geq \tfrac{\lambda_2^2}{2N}$,
\begin{eqnarray*}
\ell_N(p)-\ell_N(q) & \geq & \ell_N(\tfrac{\lambda_2^2}{2N})-\ell_N(\tfrac{\lambda_2^2}{N}) \cr
& \sim & (1-\tfrac{1}{\sqrt{2}}) \tfrac{\lambda_2}{\sqrt{N (\frac{\lambda_2^2}{2N}) \log N}} \cr
& > & \tfrac{\lambda_2}{4 \sqrt{N \xi_N \log N}}.
\end{eqnarray*}

\section{Proof of Theorem \ref{thm2}}

{\sc Proof of Theorem} \ref{thm2}(a)i.
Let $h = \lfloor \tfrac{4 (1-\epsilon) \log T}{\Delta^2 V} \rfloor$ for some $0 < \epsilon < 1$.
Let $P_0$ denote probability with respect to $\mu^n_t=0$ for all $n$ and $t$. 
Let $t_k=(2k-1)h$ and let $P_k$, 
$1 \leq k \leq K := \lfloor \tfrac{T}{2h} \rfloor$,
denote probability under which,
for $n \leq V$,
\begin{eqnarray} \label{EX}
\mu^n_{t_k-h+1} & = & \cdots = \mu^n_{t_k} = -\tfrac{\Delta}{2}, \\ \nonumber
\mu^n_{t_k+1} & = & \cdots = \mu^n_{t_k+h} = \tfrac{\Delta}{2}, \\ \nonumber
\mu^n_t & = & 0 \mbox{ for } t \leq t_k-h \mbox{ and } t > t_k+h,
\end{eqnarray}
and $\mu^n_1 = \cdots = \mu^n_T=0$ for $n>V$.
Let $E_k$ denote expectation with respect to $P_k$.  

Let $P_* = \tfrac{1}{K} \sum_{k=1}^K P_k$ and let $L = \tfrac{1}{K} \sum_{k=1}^K L_k$,
where $L_k = \tfrac{dP_k}{dP_0}(\bX)$ with $\bX = (X^n_t: 1 \leq n \leq N, 1 \leq t \leq T)$.
Hence
\begin{equation} \label{logLk}
\log L_1 = \tfrac{h \Delta}{2} \sum_{n=1}^V (\bar X^n_{h,2h}-\bar X^n_{0h})-\tfrac{h V \Delta^2}{4}.
\end{equation}

Let $A_i = \{ L \leq 3 \} \cap \{ \mbox{conclude } H_i \}$.
Since $P(A_1) = E_0(L {\bf 1}_{A_1}) \leq 3 P_0(A_1)$,
\begin{eqnarray} \label{type12}
& & \sup_{\bmu \in \Omega_0} P_{\bmu}(\mbox{Type I error}) \\ \nonumber
& & \quad + \sup_{\bmu \in \Omega_1(\Delta,V,h)} P_{\bmu}(\mbox{Type II error}) \\ \nonumber
& \geq & P_0(\mbox{conclude } H_1) + P_*(\mbox{conclude } H_0) \\ \nonumber
& \geq & P_0(A_1)+P_*(A_0) \geq \tfrac{1}{3} P_*(L \leq 3)=\tfrac{1}{3} P_1(L \leq 3),
\end{eqnarray}
with the last equality due to $L$ having the same distribution under all $P_k$ and $P_*$.

Since $E_1 L_k=1$ for $k \geq 2$,
it follows that
$P_1( \tfrac{1}{K} \sum_{k=2}^K L_k \leq 2) \geq \tfrac{1}{2}$.
Hence by (\ref{type12}),
to show that $\sup_{\bmu \in \Omega_0} P_{\bmu}(\mbox{Type I error})
+\sup_{\bmu \in \Omega_1(\Delta,V,h)} P_{\bmu}(\mbox{Type II error}) \rightarrow 0$ is not possible,
it suffices to show that
\begin{equation} \label{liminf}
P_1(L_1 \leq K) \rightarrow 1 \mbox{ as } T \rightarrow \infty. 
\end{equation}

By (\ref{logLk}),
$\log L_1 \sim {\rm N}(\tfrac{hV \Delta^2}{4},\tfrac{hV \Delta^2}{2})$,
and indeed
\begin{equation} \label{L1k}
P_1(L_1 \leq K) = \Phi \Big( \tfrac{\log K-\frac{1}{4} hV \Delta^2}{\sqrt{\frac{1}{2} hV \Delta^2}} \Big) \rightarrow 1. 
\end{equation}
$\wbox$

\medskip
{\sc Proof of Theorem} \ref{thm2}(a)ii. 
Proceed as in the proof of Theorem \ref{thm2}(a)i.,
but with $h=\lfloor \tfrac{4(1-\epsilon) \rho_Z(\beta,\zeta) \log N}{\Delta^2} \rfloor$,
and $P_k$  probability under which,
independently for $1 \leq n \leq N$,
$Q^n=1$ with probability $2N^{-\beta}$ and $Q^n=0$ otherwise.
When $Q^n=1$,
(\ref{EX}) holds.
When $Q^n=0$,
$\mu^n_1 = \cdots = \mu^n_T=0$. 

By the law of large numbers,
$P_1(\bmu \in \Omega_1(h,\Delta,V)) = P_1( \sum_{n=1}^N Q^n \geq V) \rightarrow 1$.
Hence by (\ref{type12}) it suffices to show (\ref{liminf}) with
\begin{eqnarray} \label{Lk1}
L_1 & = & \prod_{n=1}^N [1+2N^{-\beta} (e^{Z^n \Delta \sqrt{\frac{h}{2}}-\frac{h \Delta^2}{4}} -1) ], \\ \label{Zn}
Z^n & = & \sqrt{\tfrac{h}{2}}(\bar X^n_{h,2h}-\bar X^n_{0h}) \sim {\rm N} \Big( Q^n \Delta \sqrt{\tfrac{h}{2}},1 \Big).
\end{eqnarray}

\medskip
{\bf Case 1}: $\tfrac{1-\zeta}{2} < \beta < \tfrac{3(1-\zeta)}{4}$.
Recall that $\rho_Z(\beta,\zeta) = \beta-\tfrac{1-\zeta}{2}$.
By (\ref{Lk1}) and (\ref{Zn}),
\begin{eqnarray*}
E_1 L_1 & = & (1+4N^{-2 \beta}[\exp(\tfrac{h \Delta^2}{2})-1])^N \cr
& \leq & \exp(4N^{1-2 \beta+2(1-\epsilon) \rho_Z(\beta,\zeta)}) \cr
& = & \exp(4N^{\zeta-2 \epsilon \rho_Z(\beta,\zeta)}).
\end{eqnarray*}
Since $\log K =\log(\lfloor \tfrac{T}{2h} \rfloor) \sim N^{\zeta}$,
it follows that $P_1(L_1 \leq K) \geq 1-K^{-1} E_1 L_1 \rightarrow 1$ and (\ref{liminf}) holds. 

\medskip
{\bf Case 2}: $\tfrac{3(1-\zeta)}{4} \leq \beta < 1-\zeta$.
Recall that $\rho_Z(\beta,\zeta) = (\sqrt{1-\zeta}-\sqrt{1-\zeta-\beta})^2$.
Express $\log L_1 = \sum_{i=0}^3 R_i$,
where
\begin{eqnarray*}
R_i & = & \sum_{n \in \Gamma_i} \log \Big( 1+2N^{-\beta} \Big[\exp \Big( Z^n \Delta \sqrt{\tfrac{h}{2}}-\tfrac{\Delta^2 h}{4} \Big)-1 \Big] \Big), \cr
\Gamma_0 & = & \{ n: Q^n=0 \}, \cr
\Gamma_1 & = & \{ n: Q^n=1, Z^n \leq \sqrt{2(1-\zeta) \log N} \}, \cr
\Gamma_2 & = & \{ n: Q^n=1, \sqrt{2(1-\zeta) \log N} < Z^n \leq 2 \sqrt{2 \log N} \}, \cr
\Gamma_3 & = & \{ n: Q^n=1, Z^n > 2 \sqrt{2 \log N} \}.
\end{eqnarray*}
We show (\ref{liminf}) by showing that
\begin{equation} \label{Ri}
P_1(R_i \geq \tfrac{1}{4} \log K) \rightarrow 0\mbox{ for } 0 \leq i \leq 3.
\end{equation}

\medskip
${\bf i=3}$:
Since $\Delta \sqrt{\tfrac{h}{2}} \leq \sqrt{2 \log N}$,
$$P_1(R_3>0) \leq 2N^{1-\beta} \Phi(-\sqrt{2 \log N}) \rightarrow 0.
$$ 

\medskip
${\bf i=2}$:
Since
\begin{eqnarray*}
& & \Delta \sqrt{\tfrac{h}{2}} \cr
& \leq & \sqrt{2(1-\zeta) \log N}-\sqrt{2(1-\zeta-\beta) \log N} - \sqrt{2 \delta \log N}
\end{eqnarray*} 
for some $\delta>0$,
it follows that 
$$\Phi \Big( \Delta \sqrt{\tfrac{h}{2}} -\sqrt{2(1-\zeta) \log N} \Big) =o(N^{\zeta+\beta-1-\delta}).
$$
Hence
\begin{eqnarray} \label{EGamma2}
& & E_1 R_2 \\ \nonumber
& \leq & E_1(\# \Gamma_2) \log(1+2N^{4-\beta}) \\ \nonumber
& \lesssim & (N^{1-\beta} \log N) \Phi \Big( \Delta \sqrt{\tfrac{h}{2}} - \sqrt{2(1-\zeta) \log N} \Big) \\ \nonumber 
& = & o(N^{\zeta-\delta} \log N),
\end{eqnarray}
and (\ref{Ri}) follows from $\log K \sim N^{\zeta}$.

\medskip
${\bf i=1}$:
Since $\log (1+x) \leq x$,
\begin{eqnarray} \label{ER1}
& & E_1 R_1 \leq 4N^{1-2 \beta} e^{-h \Delta^2/4} \\ \nonumber
& & \qquad \times \int_{-\infty}^{\sqrt{2(1-\zeta) \log N}} \tfrac{1}{\sqrt{2 \pi}} 
e^{-(z-\Delta \sqrt{\frac{h}{2}})^2/2+z \Delta \sqrt{\frac{h}{2}}} dz \\ \nonumber
& & \qquad = 4N^{1-2\beta} \Phi \Big( \sqrt{2(1-\zeta) \log N}-2 \Delta \sqrt{\tfrac{h}{2}} \Big) e^{h \Delta^2/2} \\ \nonumber
& & \qquad 
\leq 4N^{1-2 \beta-(\sqrt{1-\zeta}-2 \sqrt{(1-\epsilon) \rho_Z(\beta,\zeta)})^2+2(1-\epsilon) \rho_Z(\beta,\zeta)} \\ \nonumber 
& & \qquad = 4N^{\zeta-\delta} \mbox{ for some } \delta>0. 
\end{eqnarray}
The last step above is shown below.
Since
$$R_1 \geq (\# \Gamma_1) \log(1-2N^{-\beta}) \stackrel{p}{\sim} -2N^{1-2 \beta} = o(N^{\zeta}),
$$
and $\log K \sim N^{\zeta}$,
(\ref{Ri}) follows from (\ref{ER1}) and Markov's inequality. 

\medskip
${\bf i=0}$:
Since $E_1 e^{R_0}=1$,
$$P_1(R_0 \geq \tfrac{1}{4} \log K) \leq K^{-\frac{1}{4}} \rightarrow 0. 
$$
$\wbox$

\medskip
{\sc Proof of} (\ref{ER1}).
It suffices to show that
\begin{equation} \label{ERprove}
1-2 \beta-(\sqrt{1-\zeta}-2 \sqrt{(1-\epsilon) \rho_Z(\beta,\zeta)})+2(1-\epsilon) \rho_Z(\beta,\zeta) < \zeta.
\end{equation}
Let $m(\rho) = -(\sqrt{1-\zeta}-2 \sqrt{\rho})^2+2 \rho$.
Inequality (\ref{ERprove}) follows from
\begin{eqnarray*}
& & m(\rho_Z(\beta,\zeta)) \cr
& = & -(\sqrt{1-\zeta}-2 \sqrt{\rho_Z(\beta,\zeta)})+2 \rho_Z(\beta,\zeta) \cr
& = & -(2 \sqrt{1-\zeta-\beta}-\sqrt{1-\zeta})^2 \cr
& & \quad + 2(\sqrt{1-\zeta}-\sqrt{1-\zeta-\beta})^2 \cr
& = & 1-\zeta-2(1-\zeta-\beta) = \zeta-1+2 \beta,
\end{eqnarray*}
and
\begin{eqnarray*}
\tfrac{d}{d \rho} m(\rho) & = & 2 \rho^{-\frac{1}{2}}(\sqrt{1-\zeta}-2 \sqrt{\rho})+2 \cr
& = & 2 \rho^{-\frac{1}{2}} \sqrt{1-\zeta}-2 > 0 \mbox{ for } \rho < 1-\zeta.
\end{eqnarray*}
$\wbox$

\medskip
{\sc Proof of Theorem} \ref{thm2}(b).
Consider $\Delta=T^{-\eta}$.
Proceed as in the proof of Theorem~\ref{thm2}(a),
with $h = \lfloor \tfrac{4(1-2 \eta)(1-\epsilon)}{V} T^{2 \eta} \log T \rfloor$ for i. and $h = \lfloor 4(1-\epsilon) \rho_Z(\beta,\zeta) T^{2 \eta} \log N \rfloor$ for ii.
For Theorem \ref{thm2}(b)i.,
\begin{equation} \label{logK}
\log K = \log(\lfloor \tfrac{T}{2h} \rfloor \sim (1-2 \eta) \log T,
\end{equation} 
and (\ref{liminf}) holds because
$$P_1(L \leq K) = \Phi \Big( \tfrac{\log K-\frac{1}{4} h V \Delta^2}{\sqrt{\frac{1}{2} hV \Delta^2}} \Big) \rightarrow 1.
$$
For Theorem \ref{thm2}(b)ii. the arguments in the proof of Theorem \ref{thm2}(a)ii. apply with
$$\log K = \log(\lfloor \tfrac{T}{2h} \rfloor) \sim (1-2 \eta) N^{\zeta}. 
$$
$\wbox$

\section{Proof of Theorem \ref{thm3}}

For $(s,t,u) \in {\cal A}_i(T)$,
the penalty of the SL scores is
$$\log(\tfrac{T}{4}(\tfrac{1}{t-s}+\tfrac{1}{u-t})) \geq \log(\tfrac{T}{2h_i}).
$$
Moreover $\#{\cal A}_i(T) \leq \tfrac{T}{d_i}$.
Hence by (\ref{Mark}) and $c_T - \log(\sum_{i=1}^{i_T} \tfrac{h_i}{d_i}) \rightarrow \infty$,
for $\bmu \in \Omega_0$,

\begin{eqnarray} \label{Pomega0}
& & P_{\bmu}(\mbox{Type I error}) \\ \nonumber
& \leq & \sum_{i=1}^{i_T} \sum_{(s,t,u) \in {\cal A}_i(T)} P_{\bmu}(\ell_N(\bp_{stu}) \geq c_T+\log(\tfrac{T}{2h_i})) \\ \nonumber
& \leq & \sum_{i=1}^{i_T} \tfrac{T}{d_i} \exp(-c_T-\log(\tfrac{T}{2h_i})) \\ \nonumber
& = & 2e^{-c_T} \sum_{i=1}^{i_T} \tfrac{h_i}{d_i} \rightarrow 0.
\end{eqnarray}

Consider $\bmu \in \Omega_1(\Delta,h,V)$ and let $\tau_j$ be the change-point satisfying the conditions in the definition of 
$\Omega_1(\Delta,h,V)$.
Let $Q^n=1$ if $|\mu_{\tau_j+1}^n-\mu_{\tau_j}^n| \geq \Delta$ and $Q^n=0$ otherwise.
We assume without loss of generality that $0 < \epsilon < 1$.

To aid in the checking of the proof of Theorem \ref{thm3},
we provide here the key ideas.
Let $j$ be such that 
$$\min(\tau_j-\tau_{j-1},\tau_{j+1}-\tau_j) \geq h \mbox{ and } m_{j \Delta} \geq V.
$$
Consider $\Delta>0$ fixed and $V \sim N^{1-\beta}$ for some $\tfrac{1-\zeta}{2} < \beta < 1-\zeta$.
Since $h \rightarrow \infty$,
it follows from (\ref{hi}) that for $N$ large 
we are able to find $(s,t,u)=(s(ik),t(ik),u(ik))$ close to $(\tau_j-h,\tau_j,\tau_j+h)$ such that
\begin{equation} \label{EZ}
E_{\bmu} Z_{stu}^n \geq [1+o(1)] \tfrac{h \Delta^2}{2} \mbox{ for } n \mbox{ satisfying } |\mu_{\tau_{j+1}}^n-\mu_{\tau_j}^n| \geq \Delta.
\end{equation}

Recall that $p_{stu}^n = 2 \Phi(-|Z_{stu}^n|)$ and let 
$q_{stu}^n = \Phi(-|Z_{stu}^n|+E_{\bmu} Z_{stu}^n)+\Phi(-|Z_{stu}^n|-E_{\bmu} Z_{stu}^n)$.
Let
\begin{eqnarray} \label{gam0}
& & \Gamma = \{ n: |Z_{stu}^n| \geq \sqrt{2 \omega \log N}, \ q_{stu}^n \geq N^{\zeta-1}, \\ \nonumber
& & \qquad \qquad |\mu_{\tau_{j+1}}^n-\mu_{\tau_j}^n| \geq \Delta \},
\end{eqnarray}
with $\omega= 1-\zeta$ when $\tfrac{3(1-\zeta)}{4} < \beta < 1-\zeta$ and $ \omega = 4(\beta-\tfrac{1-\zeta}{2})$ when
$\tfrac{1-\zeta}{2} < \beta \leq \tfrac{3(1-\zeta)}{4}$.
It follows from Lemmas \ref{lem1} and \ref{lem2} that with probability tending to 1,
\begin{eqnarray*}
\ell_N(\bp_{stu}) & \geq & \ell_N(\bq_{stu}) + (\# \Gamma) \tfrac{\lambda_2}{4 \sqrt{N \xi_N \log N}} \cr
& \geq & -\lambda_2^2 \sqrt{\log N} + (\# \Gamma) \tfrac{\lambda_2}{4 \sqrt{N \xi_N \log N}} 
\end{eqnarray*}
for $\xi_N = N^{-\omega}$.

Since the penalty $\log(\tfrac{T}{4}(\tfrac{1}{t-s}+\tfrac{1}{u-t})) \leq \log T \sim N^{\zeta}$, 
$c_T = o(\log T)$ and $\lambda_2 \sim \tfrac{N^{\frac{\zeta}{2}}}{\sqrt{\zeta \log N}}$,
to show $P_{\bmu}(\ell^{\rm pen}_{stu}(\bp) \geq c_T) \rightarrow 1$,
it suffices to show that there exists $\delta>0$ such that
\begin{equation} \label{Egamma}
E_{\bmu}(\# \Gamma) \gtrsim \left\{ \begin{array}{ll} N^{\zeta+\delta} & \mbox{ if } \tfrac{3(1-\zeta)}{4} < \beta < 1-\zeta, \cr
N^{\frac{3}{2}-2 \beta-\tfrac{\zeta}{2}+\delta} & \mbox{ if } \tfrac{1-\zeta}{2} < \beta \leq \tfrac{3(1-\zeta)}{4}.
\end{array} \right.
\end{equation}
$\wbox$

\medskip 
{\sc Proof of Theorem} \ref{thm3}(a)i. 
Consider $V =o(\tfrac{\log T}{\log N})$. 
Since $h = 4(1+\epsilon)(\tfrac{\log T}{\Delta^2 V}) \rightarrow \infty$,
$\tfrac{h_{i+1}}{h_i} \rightarrow 1$ and $d_i=o(h_i)$,
for large $T$ there exists 
$$h_i \geq 4(1+\epsilon)^{\frac{1}{2}}(\tfrac{\log T}{\Delta^2 V})
$$
such that for all $\bmu \in \Omega_1(h,\Delta,V)$, 
there exists $k$ satisfying
\begin{equation} \label{suik}
\tau_{j-1} < s(ik) < u(ik) < \tau_{j+1} \mbox{ and } |t(ik)-\tau_j| \leq \tfrac{d_i}{2}.
\end{equation}
Hence when $Q^n=1$,
\begin{equation} \label{EQn}
|E_{\bmu} Z^n_{stu}| \geq \Delta (1-\tfrac{d_i}{2h_i}) \sqrt{\tfrac{h_i}{2}} \geq \sqrt{2(1+\epsilon)^{\frac{1}{3}} V^{-1} \log T},
\end{equation}
where $(s,t,u)=(s(ik), t(ik), u(ik))$.

Let $\Gamma = \{ n: Q^n=1, |Z^n_{stu}| \geq \sqrt{2(1+\epsilon)^{\frac{1}{4}}(\tfrac{\log T}{V})} \}$.
Let $p^n_{stu} = 2 \Phi(-|Z_{stu}^n|)$ and $q_{stu}^n = \Phi(-|Z_{stu}^n|+E_{\bmu} Z_{stu}^n)+\Phi(-|Z_{stu}^n|-E_{\bmu} Z_{stu}^n)$.
Since $q^n \stackrel{\rm i.i.d.}{\sim} \mbox{Uniform}$(0,1) and $\ell_N(1) \geq -1$ for $N$ large,
by Lemmas \ref{lem1} and \ref{lem2}, 
with probability tending to 1,
\begin{eqnarray} \label{3ai}
\ell_N(\bp_{stu}) & \geq & \ell_N(\bq_{stu}) \\ \nonumber
& & \quad +(\# \Gamma) \Big[ \ell_N \Big(2 \Phi \Big( 
-\sqrt{2(1+\epsilon)^{\frac{1}{4}} \tfrac{\log T}{V}}
\Big) \Big) -1 \Big] 
\\ \nonumber
& \geq & -\lambda_2^2 \sqrt{\log N}+V[(1+\epsilon)^{\frac{1}{5}} \tfrac{\log T}{V}-\log N] \\ \nonumber
& \geq & (1+\epsilon)^{\frac{1}{6}} \log T.
\end{eqnarray}
Since the penalty $\log(\tfrac{T}{4}(\tfrac{1}{t-s}+\tfrac{1}{u-t})) \leq \log T$ and $c_T = o(\log T)$, 
it follows that $P_{\bmu}(\ell^{\rm pen}_N(\bp_{stu}) \geq c_T) \rightarrow 1$. 
$\wbox$

\medskip
{\sc Proof of Theorem} \ref{thm3}(a)ii. 
{\bf Case 1}: $V \sim N^{1-\beta}$ for $\tfrac{3(1-\zeta)}{4} \leq \beta < 1-\zeta$.
Since $h \Delta^2 = 4(1+\epsilon) (\sqrt{1-\zeta}-\sqrt{1-\zeta-\beta})^2 \log N$ and $d_i = o(h_i)$, 
for large $N$ there exists $i$ satisfying $h_i \geq (1+\epsilon)^{-\frac{1}{2}} h$ such that whenever $Q^n=1$,
\begin{eqnarray} \label{E1log}
|E_{\bmu} Z^n | & \geq & \Delta (1-\tfrac{d_i}{2h_i}) \sqrt{\tfrac{h_i}{2}} \geq \sqrt{2 \nu \log N}, \\ \nonumber
\nu & = & (1+\epsilon)^{\frac{1}{3}} (\sqrt{1-\zeta}-\sqrt{1-\zeta-\beta})^2,
\end{eqnarray}
with $(s,t,u)=(s(ik),t(ik),u(ik))$ for $k$ satisfying (\ref{suik}).

For $\Gamma$ defined in (\ref{gam0}),
\begin{eqnarray*}
E_{\bmu}(\# \Gamma) & \geq & V [\Phi \Big( -\sqrt{2(1-\zeta) \log N}+\sqrt{2 \nu \log N} \Big) \cr
& & \quad -N^{\zeta-1}] \cr
& \gtrsim & N^{1-\beta-(\sqrt{1-\zeta}-\sqrt{\nu})^2} (\log N)^{-\frac{1}{2}}, 
\end{eqnarray*}
and (\ref{Egamma}) follows from 
$$\sqrt{1-\zeta}> \sqrt{\nu} > \sqrt{1-\zeta}-\sqrt{1-\zeta-\beta}.
$$

\medskip
{\bf Case 2}: $V \sim N^{1-\beta}$ for $\tfrac{1-\zeta}{2} < \beta < \tfrac{3(1-\zeta)}{4}$.
Since $h \Delta^2 = 4(1+\epsilon) (\beta-\tfrac{1-\zeta}{2}) \log N$,
for large $N$ there exists $h_i \geq (1+\epsilon)^{-\frac{1}{2}} h$ such that whenever $Q^n=1$,
\begin{eqnarray} \label{E2}
|E_{\bmu} Z^n_{stu} | & \geq & \Delta (1-\tfrac{d_i}{2h_i}) \sqrt{\tfrac{h_i}{2}} \geq \sqrt{2 \nu \log N}, \\ \nonumber
\nu & = & (1+\epsilon)^{\frac{1}{3}} (\beta-\tfrac{1-\zeta}{2}),
\end{eqnarray}
with $(s,t,u) = (s(ik),t(ik),u(ik))$ for $k$ satisfying (\ref{suik}).

For $\Gamma$ defined in (\ref{gam0}),
\begin{eqnarray*}
E_{\bmu}(\# \Gamma) & \geq & V [\Phi \Big( -2 \sqrt{(2 \beta-1+\zeta) \log N} + \sqrt{2 \nu \log N} \Big) \cr
& & \quad -N^{\zeta-1}] \cr
& \gtrsim & N^{1-\beta-(2 \sqrt{\beta-\frac{1-\zeta}{2}}-\sqrt{\nu})^2} (\log N)^{-\frac{1}{2}}, 
\end{eqnarray*}
and (\ref{Egamma}) follows from 
$$2 \sqrt{\beta-\tfrac{1-\zeta}{2}} > \sqrt{\nu} > \sqrt{\beta-\tfrac{1-\zeta}{2}}. 
$$ 
$\wbox$

\medskip
{\sc Proof of Theorem} \ref{thm3}(b).
Consider first $V = o(\tfrac{\log T}{\log N})$.
Let $h_i \geq (1-2 \eta)(1+\epsilon)^{\frac{1}{2}} (\tfrac{\log T}{\Delta^2 V})$ be such that for all $\bmu \in \Omega_1(h,\Delta,V)$,
(\ref{suik}) holds for some $k$.
Let
$$\Gamma = \{ n: Q^n=1, |Z^n_{stu}| \geq \sqrt{2(1-2 \eta)(1+\epsilon)^{\frac{1}{4}} \tfrac{\log T}{V}} \}
$$
and define $p_{stu}^n$ and $q_{stu}^n$ as in the proof of Theorem \ref{thm3}(a)i.

By the arguments in (\ref{3ai}),
with probability tending to 1,
$$\ell_N(\bp_{stu}) \geq (1-2 \eta)(1+\epsilon)^{\frac{1}{6}} \log T.
$$
Since $\Delta \sim T^{-\eta}$,
it follows that $h_i \gtrsim T^{2 \eta} \log N$ and the penalty 
$\log(\tfrac{T}{4}(\tfrac{1}{u-t}+\tfrac{1}{t-s})) \leq (1-2 \eta) \log T$ for $T$ large.
Hence by $c_T=o(\log T)$ we conclude $P_{\bmu}(\ell_N^{\rm pen}(\bp_{stu}) \geq c_T) \rightarrow 1$.

For $V \sim N^{1-\beta}$ with $\tfrac{1-\zeta}{2} < \beta < 1-\zeta$,
the same arguments in the proof of Theorem~\ref{thm3}(a)ii. apply.  
$\wbox$

\section{Proof of Theorem \ref{thm4}}

{\sc Proof of Theorem} \ref{thm4}(a).
Let $h = \lfloor \tfrac{(1-\epsilon) \log T}{\mu_0 V I_r} \rfloor$ for some $0 < \epsilon < 1$.
Let $P_0$ denote probability with respect to $\mu_t^n = (\tfrac{1+r}{2}) \mu_0$ for all $n$ and $t$.
Let $t_k = (2k-1)h$.
Let $P_k$,
$1 \leq k \leq K:=\lfloor \tfrac{T}{2h} \rfloor$,
denote probability under which for $n \leq V$,
\begin{equation} \label{lamn}
\mu^n_t = \left\{ \begin{array}{ll} \mu_0 & \mbox{ for } t_k-h < t \leq t_k, \cr
r \mu_0 & \mbox{ for } t_k < t \leq t_k+h, \cr
(\tfrac{1+r}{2}) \mu_0 & \mbox{ for } t \leq t_k-h \mbox{ and } t > t_k+h, \end{array} \right.
\end{equation}
and $\mu^n_1 = \cdots = \mu^n_T = (\tfrac{1+r}{2}) \mu_0$ for $n>V$. 
Let $E_k$ and ${\rm Var}_k$ denote expectation and variance respectively with respect to $P_k$.
Let
\begin{eqnarray} \label{Ukn}
U^n & = & S_{0h}^n \log(\tfrac{2}{1+r}) +S_{h,2h}^n \log(\tfrac{2r}{r+1}), \\ \label{Lk}
L_1 & = & \tfrac{dP_1}{dP_0}(\bX) = \prod_{n=1}^V \exp(U^n).
\end{eqnarray}

By (\ref{type12})--(\ref{liminf}),
it suffices to show that
\begin{equation} \label{lim0}
P_1(L_1 \leq K) \rightarrow 1 \mbox{ as } T \rightarrow \infty.
\end{equation}
Since $E_1 (\log L_1) = h \mu_0 VI_r$ and ${\rm Var}_1(\log L_1) = h \mu_0 V C_r$, 
where $C_r = r [\log(\tfrac{2r}{r+1})]^2+[\log(\tfrac{2}{r+1})]^2$,
by Chebyshev's inequality,
$$P_1(L_1 \leq K) \geq 1-\tfrac{hV \mu_0 C_r}{(\log K-hV \mu_0 I_r)^2} \rightarrow 1,
$$
and (\ref{lim0}) holds. 
$\wbox$

\medskip
We preface the proof of Theorem \ref{thm4}(b) with Lemma \ref{lem3},
which provides an alternative representation of $\rho_r(\beta,\zeta)$.
Let
\begin{eqnarray} \label{Domega}
D(\omega) & = & \tfrac{1}{1+r^{\omega}} \log(\tfrac{2}{1+r^{\omega}})+\tfrac{r^{\omega}}{1+r^{\omega}} \log(\tfrac{2r^{\omega}}{1+r^{\omega}}), \\ \nonumber
g(\omega) & = & (\tfrac{1+r^{\omega}}{2})^{\frac{1}{\omega}}.
\end{eqnarray}
Let $\xi(\omega) = \tfrac{\beta-\omega^{-1}(1-\zeta)}{2 g(\omega)-1-r}$.
Recall from (\ref{rho}) that 
\begin{equation} \label{rholem}
\rho_r(\beta,\zeta) = \max_{\frac{1-\zeta}{\beta} < \omega \leq 2} \xi(\omega) \mbox{ for } \tfrac{1-\zeta}{2} < \beta < 1-\zeta. 
\end{equation}

\medskip
\begin{lem} \label{lem3}
For $\tfrac{1}{2} < \tfrac{\beta}{1-\zeta} \leq \tfrac{1}{2}[1+\tfrac{2g(2)-1-r}{g(2) D(2)}]$,
$\xi$ achieves its maximum at $\omega=2$ and
\begin{equation} \label{rho3}
\rho_r(\beta,\zeta) = \tfrac{\beta-\frac{1}{2}(1-\zeta)}{2 g(2)-1-r}.
\end{equation}
For $\tfrac{1}{2}[1+\tfrac{2g(2)-1-r}{g(2) D(2)}] < \tfrac{\beta}{1-\zeta} < 1$,
$\xi$ achieves its maximum at some $\omega<2$ and 
\begin{equation} \label{rho2}
\rho_r(\beta,\zeta) = \tfrac{1-\zeta}{2 g(\omega) D(\omega)}. 
\end{equation}
\end{lem}

\medskip
{\sc Proof}.
Since
\begin{eqnarray*}
\tfrac{d}{d \omega} \log \xi(\omega) & = & \tfrac{\omega^{-2}(1-\zeta)}{\beta-\omega^{-1}(1-\zeta)} - \tfrac{2 \frac{d}{d \omega} g(\omega)}{2 g(\omega)-1-r}, \cr
\tfrac{d}{d \omega} g(\omega) & = & \tfrac{d}{d \omega} \exp[\tfrac{1}{\omega} \log(\tfrac{1+r^{\omega}}{2})] \cr
& = & [\tfrac{r^{\omega} \log r}{\omega(1+r^{\omega})}-\tfrac{1}{\omega^2} \log(\tfrac{1+r^{\omega}}{2})] g(\omega) \cr
& = & \tfrac{D(\omega) g(\omega)}{\omega^2},
\end{eqnarray*}
it follows that $\tfrac{d}{d \omega} \log \xi(\omega)=0$ when 
\begin{equation} \label{betaomega}
\omega^{-2}(1-\zeta)[2 g(\omega)-1-r] = 2[\beta-\omega^{-1}(1-\zeta)] \tfrac{D(\omega) g(\omega)}{\omega^2},
\end{equation}
that is when 
\begin{equation} \label{D11}
\tfrac{\beta}{1-\zeta} = \omega^{-1} + \tfrac{2 g(\omega)-1-r}{2 g(\omega) D(\omega)}.
\end{equation}
For $\tfrac{1}{2} < \tfrac{\beta}{1-\zeta} \leq \tfrac{1}{2}[1+\tfrac{2 g(2)-1-r}{g(2) D(2)}]$,
the solution of $\omega$ to (\ref{D11}) is at least 2 and the maximum in (\ref{rholem}) is attained at $\omega=2$.
We conclude (\ref{rho3}).
For $\tfrac{1}{2}[1+\tfrac{2g(2)-1-r}{g(2) D(2)}] < \tfrac{\beta}{1-\zeta} < 1$,
the solution of $\omega$ to (\ref{D11}) lies in the interval $(\frac{1-\zeta}{\beta},2)$.
We conclude (\ref{rho2}) from (\ref{rholem}) and a rearrangement of (\ref{betaomega}).
$\wbox$

\medskip
{\sc Proof of Theorem} \ref{thm4}(b).
For $\frac{1-\zeta}{2} < \beta < 1-\zeta$,
let $\omega$ be the maximizer in 
\begin{equation} \label{rhoY}
\rho_r(\beta,\zeta) = \max_{\frac{1-\zeta}{\beta} < \omega \leq 2} (\tfrac{\beta-\omega^{-1}(1-\zeta)}{2 g(\omega)-1-r}).
\end{equation}

Let $h= \lfloor \tfrac{(1-\epsilon) \rho_r(\beta,\zeta) \log N}{\mu_0} \rfloor$ for some $\epsilon>0$.
Let $P_0$ denote probability with respect to $\mu^n_t = g(\omega) \mu_0$ for all $n$ and $t$.
Let $t_k = (2k-1)h$.
Let $P_k$,
$1 \leq k \leq K:= \lfloor \tfrac{T}{2h} \rfloor$,
denote probability under which,
independently for $1 \leq n \leq N$,
$Q^n=1$ with probability $2N^{-\beta}$,
and $Q^n=0$ otherwise.
When $Q^n=1$,
\begin{equation} \label{mutn}
\mu^n_t = \left\{ \begin{array}{ll} \mu_0 & \mbox{ for } t_k-h < t \leq t_k, \cr
r \mu_0 & \mbox{ for } t_k < t \leq t_k+h, \cr
g(\omega) \mu_0 & \mbox{ for } t \leq t_k-h \mbox{ and } t > t_k+h. \end{array} \right.
\end{equation}
When $Q^n=0$,
$\mu_1^n = \cdots  =\mu_T^n =g(\omega) \mu_0$.
Let $E_1$ denote expectation with respect to $P_1$.
Let $P_Q = P_1(\cdot|Q^1=1)$ and let $E_Q$ denote expectation with respect to $P_Q$.  

By (\ref{type12})--(\ref{liminf}),
it suffices to show (\ref{lim0}) for
\begin{eqnarray} \label{L14}
L_1 & = & \tfrac{dP_1}{dP_0}(\bX) = \prod_{n=1}^N (1+2N^{-\beta}[\exp(U^n)-1]), \\ \nonumber
U^n & = & S^n_{0h} \log(\tfrac{1}{g(\omega)})+S^n_{h,2h} \log(\tfrac{r}{g(\omega)}) \\ \nonumber
& & \quad -h \mu_0[1+r-2 g(\omega)].
\end{eqnarray}
For notational simplicity,
let $S_{0h} = S_{0h}^1$ and $S_{h,2h} = S_{h,2h}^1$.

For $X \sim$Poisson($\lambda$) and constant $C>0$,
\begin{equation} \label{Poi}
E(C^X) = \sum_{x=0}^{\infty} e^{-\lambda} \tfrac{(C \lambda)^x}{x!} = e^{\lambda(C-1)}.
\end{equation}
This identity is applied in (\ref{EUn}),
(\ref{E11e}) and (\ref{E1G}).

\medskip
{\bf Case 1}: $\frac{1}{2} < \tfrac{\beta}{1-\zeta} \leq \frac{1}{2}[1+\tfrac{2 g(2)-1-r}{g(2) D(2)}]$,
$\omega=2$.
By Lemma \ref{lem2},
(\ref{mutn})--(\ref{Poi}) and $[g(2)]^2=\tfrac{1+r^2}{2}$,
\begin{eqnarray} \label{EUn}
& & E_Q \exp(U^1) \\ \nonumber
& = & E_Q[(\tfrac{1}{g(2)})^{S_{0h}} (\tfrac{r}{g(2)})^{S_{h,2h}}] e^{-h \mu_0[1+r-2 g(2)]} \\ \nonumber
& = & \exp (h \mu_0[\tfrac{1}{g(2)}-1+\tfrac{r^2}{g(2)}-r]-h \mu_0[1+r-2 g(2)]) \\ \nonumber
& = & \exp (2h \mu_0[2g(2)-1-r]) \\ \nonumber
& = & \exp(\tfrac{h \mu_0 (2 \beta-1+\zeta)}{\rho_r(\beta,\zeta)}) \leq N^{(1-\epsilon)(2 \beta-1+\zeta)}.
\end{eqnarray}
Hence 
\begin{eqnarray*}
E_1 L_1 & = & (1+4N^{-2 \beta} [E_Q \exp(U^1)-1])^N \cr
& \leq & \exp(4N^{\zeta-\delta}) = o(K),
\end{eqnarray*}
where $\delta = \epsilon(2 \beta-1+\zeta)$, 
and (\ref{lim0}) holds.

\medskip
{\bf Case 2}: $\tfrac{1}{2}[1+\tfrac{2 g(2)-1-r}{g(2) D(2)}] < \tfrac{\beta}{1-\zeta} < 1$.
Express
\begin{eqnarray} \label{logL1}
\log L_1 & = & R_0 + R_1, \\ \nonumber 
\mbox{where } R_i & = & \sum_{n \in \Gamma_i} \log(1+2N^{-\beta}[\exp(U^n)-1]), \\ \nonumber  
\Gamma_0 & = & \{ n: Q^n=0 \} \cup \{ n: Q^n=1, \exp(U^n) \leq N^{\beta} \}, \\ \nonumber
\Gamma_1 & = & \{ n: Q^n=1, \exp(U^n) > N^{\beta} \}.
\end{eqnarray}
We conclude (\ref{lim0}) from 
\begin{equation} \label{check1}
P_1(R_i \leq \tfrac{1}{2} \log K) \rightarrow 1 \mbox{ for } i=0 \mbox{ and } 1.
\end{equation}

\medskip ${\bf i=0}$:
Let $a=\omega-1$ with $\omega$ the maximizer in (\ref{rhoY}).
Since $g(\omega) = \tfrac{1+r^{a+1}}{2 g^a(\omega)}$,
by (\ref{rhoY}), (\ref{L14}) and (\ref{Poi}), 
\begin{eqnarray} \label{E11e}
& & E_Q [\exp(U^1) {\bf 1}_{\{ 1 \in \Gamma_0 \}}] \\ \nonumber
& \leq & N^{\beta(1-a)} E_Q \exp(aU^1) \\ \nonumber 
& = & N^{\beta(1-a)} \exp(h \mu_0[\tfrac{1}{g^a(\omega)}-1+\tfrac{r^{a+1}}{g^a(\omega)}-r] \\ \nonumber
& & -a h \mu_0[1+r-2 g(\omega)]) \\ \nonumber
& = & N^{\beta(1-a)} \exp(\omega h \mu_0[2 g(\omega)-1-r]) \\ \nonumber
& = & N^{\beta(1-a)} \exp(\tfrac{h \mu_0 (\beta \omega-1+\zeta)}{\rho_r(\beta,\zeta)}) \leq N^{2 \beta-1+\zeta-\delta},
\end{eqnarray}
where $\delta = \epsilon (\beta \omega-1+\zeta)$.
Since $E_0 \exp(U^n)=1$,
it follows from (\ref{E11e}) that
\begin{eqnarray*}
E_1 \exp(R_0) & \leq & (1+4N^{-2 \beta} E_Q[\exp(U^1) {\bf 1}_{ \{ 1 \in \Gamma_0 \}}] )^N \cr
& \leq & \exp(4N^{\zeta-\delta}),
\end{eqnarray*}
and (\ref{check1}) holds.

\medskip
${\bf i=1}$: 
Express $U^1 = v_1 S_{0h}+v_2 S_{h,2h}-z$,
where $v_1 = \log(\tfrac{1}{g(\omega)})$, 
$v_2 = \log(\tfrac{r}{g(\omega)})$ and $z=h \mu_0[1+r-2g(\omega)]$.
Since $g(\omega) = \tfrac{1+r^{a+1}}{2 g^a(\omega)}$,
by Markov's inequality and (\ref{Poi}),
\begin{eqnarray} \label{E1G}
& & E_1(\# \Gamma_1) \\ \nonumber
& = & 2N^{1-\beta} P_Q (e^{aU^1} > N^{a \beta}) \\ \nonumber
& \leq & 2N^{1-\beta -a \beta} e^{-az} E_Q(e^{v_1 a S_{0h}} e^{v_2 a S_{h,2h}}) \\ \nonumber
& = & 2N^{1-\omega \beta} \exp(-az+h \mu_0[e^{v_1 a}-1+re^{v_2 a}-r]) \\ \nonumber
& = & 2N^{1-\omega \beta} \exp(\omega h \mu_0[2 g(\omega)-1-r]) \\ \nonumber
& = & 2N^{1-\omega \beta} \exp(\tfrac{h \mu_0 (\beta \omega-1+\zeta)}{\rho_r(\beta,\zeta)}) \leq N^{\zeta-\delta},
\end{eqnarray}
where $\delta = \epsilon (\beta \omega-1+\zeta)$.
Since 
$$R_1 \leq (\# \Gamma_1) \max_{n \in \Gamma_1} U_1^n \mbox{ and } P_1(\max_n U^n \geq N^{\frac{\delta}{2}}) \rightarrow 0,
$$
we conclude (\ref{check1}) from (\ref{E1G}) and Markov's inequality.
$\wbox$

\section{Proof of Theorem \ref{thm5}}

It follows from (\ref{Pomega0}) that
$\sup_{\bmu \in \Lambda_0} P_{\bmu}(\mbox{Type I error}) \rightarrow 0$.

Consider $\bmu \in \Lambda_1(h,\Delta,V)$ and let $\tau_j$ be a change-point such that 
$$\min(\tau_{j+1}-\tau_j,\tau_j-\tau_{j-1}) \geq h \mbox{ and } m_{j \Delta} \geq V,
$$
where $m_{j \Delta} = \# \{ n: |\log(\mu_{\tau_j+1}^n/\mu_{\tau_j}^n)| \geq \Delta \}$.

Let $Q^n=1$ if $|\log(\mu_{\tau_j+1}^n/\mu_{\tau_j}^n)| \geq \Delta$ and $Q^n=0$ otherwise.

\medskip
{\sc Proof of Theorem} \ref{thm5}(a).
Consider $V=o(\tfrac{\log T}{\log N})$ and recall from (\ref{IY}) that $I_r = r \log(\tfrac{2r}{r+1})+\log(\tfrac{2}{r+1})$.
Let $r_1$ and $\mu_1$ be such that $e^{\Delta} > r_1 > r$ and $\mu_0/(1+\epsilon)^{\frac{1}{3}} < \mu_1 < \mu_0$.
Since $h V I_r \mu_0 = (1+\epsilon) \log T$,
$\frac{h_{i+1}}{h_i} \rightarrow 1$ and $d_i=o(h_i)$,
for $T$ large there exists 
\begin{equation} \label{hi5}
h_i \geq (1+\epsilon)^{\frac{1}{2}} I_r^{-1} (\tfrac{\log T}{\mu_1 V}),
\end{equation}
such that for all $\bmu \in \Lambda_1(h,\Delta,V)$, 
there exists $k$ such that
$$\tau_{j-1} < s(ik) < u(ik) < \tau_{j+1}, \qquad |t(ik)-\tau_j| \leq \tfrac{d_i}{2}.
$$
Moreover when $Q^n=1$,
\begin{equation} \label{logr1}
|\log(E_{\bmu} Y^n_{tu}/E_{\bmu} Y^n_{st})| \geq \log r_1,
\end{equation}
where $(s,t,u) = (s(ik),t(ik),u(ik))$.
Let 
$$\Gamma = \{ n: Q^n=1, \ Y^n_{su} \geq (1+r) h_i \mu_1, \ | \log (Y_{tu}^n/Y_{st}^n) | \geq \log r \}.
$$

By (\ref{hi5}), 
for $n \in \Gamma$, 
\begin{equation} \label{pnV}
p^n_{stu} \leq 2 \exp(-\mu_1 h_i I_r) \leq 2 \exp(-(1+\epsilon)^{\frac{1}{2}} \tfrac{\log T}{V}).
\end{equation}
Since $\tfrac{\log T}{V \log N} \rightarrow \infty$,
for $N$ large,
$$\ell(p_{stu}^n) \geq (1+\epsilon)^{\frac{1}{3}}(\tfrac{\log T}{V}).
$$
Hence as $\ell_N(q) \geq -1$ for $N$ large,
by Lemma \ref{lem1},
with probability tending to 1,
\begin{eqnarray*}
\ell_N(\bp_{stu}) & \geq & \ell_N(\bq_{stu}) +(\# \Gamma)[(1+\epsilon)^{\frac{1}{3}}(\tfrac{\log T}{V})-1] \cr
& \geq & -\lambda_2^2 \sqrt{\log N}+(1+\epsilon)^{\frac{1}{4}} \log T.
\end{eqnarray*}
Since the penalty $\log(\tfrac{T}{4}(\tfrac{1}{t-s}+\tfrac{1}{u-t})) \leq \log T$,
$\lambda_2^2 \sqrt{\log N}=o(\log T)$ and $c_T=o(\log T)$,
we can conclude that
$P_{\bmu}(\ell^{\rm pen}_N(\bp_{stu}) \geq c_T) \rightarrow 1$. 
$\wbox$

\medskip
{\sc Proof of Theorem} \ref{thm5}(b).
Consider $V \sim N^{1-\beta}$ for $\tfrac{1-\zeta}{2} < \beta < 1-\zeta$.
For $N$ large,
there exists 
\begin{equation} \label{hi2}
\log N \gtrsim h_i \geq (1+\epsilon)^{\frac{1}{2}} \rho_r(\beta,\zeta) (\tfrac{\log N}{\mu_0}) 
\end{equation}
such that for all $\bmu \in \Lambda_1(h,\Delta,V)$, 
there exists $k$ such that
$$\tau_{j-1} < s(ik) < u(ik) < \tau_{j+1}, \quad |t(ik)-\tau_j| \leq \tfrac{d_i}{2},
$$ 
and conditioned on $Q^n=1$,
either
\begin{equation} \label{Q4}
E_{\bmu} Y_{tu}^n \geq r E_{\bmu} Y_{st}^n \mbox{ or } E_{\bmu} Y_{st}^n \geq r E_{\bmu} Y_{tu}^n,
\end{equation}
where $(s,t,u)=(s(ik),t(ik),u(ik))$.

By Stirling's approximation $x! \sim \sqrt{2 \pi x}(\tfrac{x}{e})^x$, 
for $X \sim$Poisson($\eta$),
as $x \rightarrow \infty$,
\begin{equation} \label{stir}
P(X=x) = e^{-\eta} \tfrac{\eta^x}{x!} \sim \tfrac{1}{\sqrt{2 \pi x}} \exp[-\eta+x-x \log(\tfrac{x}{\eta})].
\end{equation}
By apply this in (\ref{Y1}) and (\ref{Y2}).

\medskip
{\bf Case 1}: $\tfrac{1}{2} < \tfrac{\beta}{1-\zeta} \leq \tfrac{1}{2}[1+\tfrac{2 g(2)-1-r}{g(2) D(2)}]$ and $\rho_r(\beta,\zeta) = \tfrac{\beta-\frac{1}{2}(1-\zeta)}{2g(2)-1-r}$.
Let
\begin{eqnarray} \label{gam51}
\Gamma & = & \{ n: Q^n=1, \ Y_{su}^n \geq \sqrt{2(1+r^2)} h_i \mu_0-1, \\ \nonumber
& & \qquad |\log(Y_{tu}^n/Y_{st}^n)| \geq 2 \log r, \ q_{stu}^n \geq N^{\zeta-1} \}.
\end{eqnarray}
Consider $Y_1 \sim$ Poisson($h_i \mu_0$) and $Y_2 \sim$ Poisson($r h_i \mu_0$).
By (\ref{stir}) and $h_i \lesssim \log N$, 
\begin{eqnarray} \label{Y1}
& & P(Y_1 = \lfloor (\tfrac{2}{1+r^2})^{\frac{1}{2}} h_i \mu_0 \rfloor) \\ \nonumber
& \gtrsim & \tfrac{1}{\sqrt{\log N}} \exp (h_i \mu_0[-1+(\tfrac{2}{1+r^2})^{\frac{1}{2}} \\ \nonumber
& & \quad -(\tfrac{2}{1+r^2})^{\frac{1}{2}} \log((\tfrac{2}{1+r^2})^{\frac{1}{2}})]), \\ \nonumber 
& & P(Y_2 = \lceil (\tfrac{2}{1+r^2})^{\frac{1}{2}} r^2 h_i \mu_0 \rceil) \\ \nonumber
& \gtrsim & \tfrac{1}{\sqrt{\log N}} \exp (h_i \mu_0[-r+r^2(\tfrac{2}{1+r^2})^{\frac{1}{2}} \\ \nonumber
& & \quad -r^2(\tfrac{2}{1+r^2})^{\frac{1}{2}} \log(r(\tfrac{2}{1+r^2})^{\frac{1}{2}})]).
\end{eqnarray}

Recall that $g(2) = (\tfrac{1+r^2}{2})^{\frac{1}{2}}$ and 
$D(2) = \tfrac{1}{1+r^2} \log(\tfrac{2}{1+r^2})+\tfrac{r^2}{1+r^2} \log(\tfrac{2r^2}{1+r^2})$ [see (\ref{Domega})].
By (\ref{Y1}),
\begin{eqnarray} \label{PC2}
& & E_{\bmu}(\# \Gamma) \\ \nonumber
& \geq & V [P(Y_1 = \lfloor (\tfrac{2}{1+r^2})^{\frac{1}{2}} h_i \mu_0 \rfloor) P(Y_2 = \lceil (\tfrac{2}{1+r^2})^{\frac{1}{2}} r^2 h_i \mu_0 \rceil) \\ \nonumber
& & \quad -N^{\zeta-1}] \\ \nonumber
& \gtrsim & \tfrac{N^{1-\beta}}{\log N} \exp (h_i \mu_0 [2g(2)-1-r-g(2) D(2)]).
\end{eqnarray}
By (\ref{gam51}), 
for $n \in \Gamma$,
\begin{eqnarray} \label{P2n}
p^n_{stu} & \leq & 2 \exp(-Y_{su}^n D(2)) \leq \xi_N \\ \nonumber
\mbox{where } \xi_N & = & C_2 \exp(-2h_i \mu_0 g(2) D(2)) \mbox{ for } C_2 = 2e^{D(2)}.
\end{eqnarray}

Let $q^n_{stu}=F_n(p^n_{stu})$ where $F_n$ is the distribution function of $p^n_{stu}$.
It follows from Lemmas \ref{lem1} and \ref{lem2} that with probability tending to 1,
\begin{eqnarray} \label{lN4}
& &  \ell_N(\bp_{stu}) \\ \nonumber
& \geq & \ell_N(\bq_{stu}) + (\# \Gamma) \tfrac{\lambda_2}{4 \sqrt{N \xi_N \log N}} \\ \nonumber
& \geq & -\lambda_2^2 \sqrt{\log N} + \tfrac{\lambda_2 N^{\frac{1}{2}-\beta}}{(\log N)^{\frac{3}{2}}}
\exp(h_i \mu_0[2 g(2)-1-r]).
\end{eqnarray}
Since $\lambda_2 \asymp \tfrac{N^{\frac{\zeta}{2}}}{\sqrt{\log N}}$ and by (\ref{hi2})
$$h_i \mu_0 \geq (1+\epsilon)^{\frac{1}{2}} \rho_r(\beta,\zeta) \log N = (1+\epsilon)^{\frac{1}{3}}
(\tfrac{\beta-\tfrac{1}{2}(1-\zeta)}{2 g(2)-1-r})\log N,
$$
it follows from (\ref{lN4}) that $\ell_N(\bp_{stu}) \gtrsim \tfrac{N^{\zeta+\delta}}{(\log N)^2}$ for $\delta=[(1+\epsilon)^{\frac{1}{2}}-1]
[\beta-\tfrac{1}{2}(1-\zeta)]$.
Since the penalty $\log(\tfrac{T}{4}(\tfrac{1}{t-s}+\tfrac{1}{u-t})) \leq \log T \sim N^{\zeta}$ and $c_T=o(N^{\zeta})$,
we conclude $P_{\bmu}(\ell_N^{\rm pen}(\bp_{stu}) \geq c_T) \rightarrow 1$.

\medskip
{\bf Case 2}: $\tfrac{1}{2}[1+\tfrac{2g(2)-1-r}{g(2)D(2)}] < \tfrac{\beta}{1-\zeta} < 1$ and 
$\rho_r(\beta,\zeta) = \tfrac{1-\zeta}{2g(\omega) D(\omega)} = \tfrac{\beta-\omega^{-1}(1-\zeta)}{2g(\omega)-1-r}$ with $\omega$ achieving the maximum in (\ref{rholem}).
Let 
\begin{eqnarray} \label{gam52}
\Gamma & = & \{ n: Q^n=1, \ Y^n_{su} \geq 2 g(\omega) h_i \mu_0-1, \\ \nonumber
& & \qquad |\log(Y^n_{tu}/Y^n_{st})| \geq \omega \log r, \ q_{stu}^n \geq N^{\zeta-1} \}.
\end{eqnarray}
Consider $Y_1 \sim$ Poisson($h_i \mu_0$) and $Y_2 \sim$ Poisson($rh_i \mu_0$).

By (\ref{stir}) and $h_i \lesssim \log N$,
\begin{eqnarray} \label{Y2}
& & P(Y_1 = \lfloor \tfrac{2g(\omega)}{r^{\omega}+1} h_i \mu_0 \rfloor) \\ \nonumber
& \gtrsim & \tfrac{1}{\sqrt{\log N}} \exp(h_i \mu_0[-1+\tfrac{2g(\omega)}{r^{\omega}+1} \\ \nonumber
& & \quad -\tfrac{2 g(\omega)}{r^{\omega}+1} \log(\tfrac{2 g(\omega)}{r^{\omega}+1})]), \\ \nonumber
& & P(Y_2 = \lceil \tfrac{2 r^{\omega} g(\omega)}{r^{\omega}+1} h_i \mu_0 \rceil) \\ \nonumber
& \gtrsim & \tfrac{1}{\sqrt{\log N}} \exp(h_i \mu_0[-r+\tfrac{2r^{\omega} g(\omega)}{r^{\omega}+1} \\ \nonumber
& & \quad -\tfrac{2 r^{\omega} g(\omega)}{r^{\omega}+1} \log(\tfrac{2 r^{\omega-1} g(\omega)}{r^{\omega}+1})]). 
\end{eqnarray}

Recall that $g(\omega) = (\tfrac{1+r^{\omega}}{2})^{\frac{1}{\omega}}$ and
$D(\omega) = \tfrac{1}{1+r^{\omega}} \log(\tfrac{2}{1+r^{\omega}})+\tfrac{r^{\omega}}{1+r^{\omega}} \log(\tfrac{2r^{\omega}}{1+r^{\omega}})$ [see (\ref{Domega})].
By (\ref{Y2}),
\begin{eqnarray} \label{EB2}
& & E_{\bmu} (\# \Gamma) \\ \nonumber
& \geq & V [P(Y_1 = \lfloor \tfrac{2g(\omega)}{r^{\omega}+1} h_i \mu_0 \rfloor) 
P(Y_2 \geq \lceil \tfrac{2 r^{\omega} g(\omega)}{r^{\omega}+1} h_i \mu_0 \rceil) \\ \nonumber
& & \quad -N^{\zeta-1}] \\ \nonumber
& \gtrsim & \tfrac{N^{1-\beta}}{\log N}  \exp(h_i \mu_0[2 g(\omega)-1-r-2(\tfrac{\omega-1}{\omega}) g(\omega) D(\omega)]).
\end{eqnarray}
By (\ref{gam52}), 
for $n \in \Gamma$,
\begin{eqnarray} \label{B2p}
p^n & \leq & 2 \exp(-Y^n_{s(ik),u(ik)} D(\omega)) \leq \xi_N \\ \nonumber
\mbox{where } \xi_N & = & C_{\omega} \exp(-2h_i \mu_0 g(\omega) D(\omega)) \mbox{ for } C_{\omega} = 2e^{D(\omega)}.
\end{eqnarray}

Let $q^n = F_n(p^n)$
where $F_n$ is the distribution function of $p^n$.
It follows from Lemmas \ref{lem1} and \ref{lem2} that with probability tending to 1,
\begin{eqnarray} \label{lN5}
& & \ell_N(\bp_{stu}) \\ \nonumber
& \geq & \ell_N(\bq_{stu}) + (\# \Gamma) \tfrac{ \lambda_2}{4 \sqrt{N \xi_N \log N}} \\ \nonumber
& \gtrsim & -\lambda_2^2 \sqrt{\log N} + \tfrac{\lambda_2 N^{\frac{1}{2}-\beta}}{(\log N)^{\frac{3}{2}}} \\ \nonumber
& & \qquad \times \exp(h_i \mu_0[2 g(\omega)-1-r-(\tfrac{\omega-2}{\omega}) g(\omega) D(\omega)]).
\end{eqnarray}

Since $\lambda_2 \asymp \tfrac{N^{\frac{\zeta}{2}}}{\sqrt{\log N}}$ and by (\ref{hi2}),
\begin{eqnarray*}
h_i \mu_0 & \geq & (1+\epsilon)^{\frac{1}{2}} \rho_r(\beta,\zeta) \log N \cr
& = & (1+\epsilon)^{\frac{1}{2}} (\tfrac{1-\zeta}{2 g(\omega) D(\omega)}) \log N \cr
& = & (1+\epsilon)^{\frac{1}{2}} (\tfrac{\beta-\omega^{-1}(1-\zeta)}{2 g(\omega)-1-r}) \log N, 
\end{eqnarray*}
it follows from (\ref{lN5}) that $\ell_N(\bp_{stu}) \gtrsim \tfrac{N^{\zeta+\delta}}{(\log N)^2}$ for 
$\delta = [(1+\epsilon)^{\frac{1}{2}}-1][\beta-\tfrac{1}{2}(1-\zeta)]$.
Since the penalty $\log(\tfrac{T}{4}(\tfrac{1}{t-s}+\tfrac{1}{u-t})) \leq \log T \sim N^{\zeta}$ and $c_T = o(\log T)$,
we conclude $P_{\bmu}(\ell_N^{\rm pen}(\bp_{stu}) \geq c_T) \rightarrow 1$.
$\wbox$
\end{appendix}


\end{document}